\documentclass[]{amsart}
\usepackage{amsmath,amsfonts,amssymb,amsthm}
\usepackage{enumerate}
\usepackage[usenames]{color}
\usepackage{graphicx}

 \newcommand{\obra}[3]{{\sc #1} {\em #2}. {#3}.}

\newtheorem{theorem}{\bf Theorem}
 \newtheorem{lemma}[theorem]{\bf Lemma}
 \newtheorem{proposition}[theorem]{\bf Proposition}
 \newtheorem{definition}[theorem]{\bf Definition}
 
 \newtheorem{remark}[theorem]{Remark}

 \renewenvironment{proof}{{\em Proof \/.-}}
   {\hfill $\square$\newline}

\DeclareMathOperator{\Sing}{Sing}
\DeclareMathOperator{\Reg}{Reg}
\DeclareMathOperator{\Sat}{Sat}

 \newcommand{\R}{\mathbb{R}}
 \newcommand{\Q}{\mathbb{Q}}
 \newcommand{\N}{\mathbb{N}}

 \newcommand{\SSS}{\mathbb{S}}

 \newcommand{\CC}{\mathcal{C}}
  
  \newcommand{\GG}{\mathcal{G}}

    \newcommand{\LL}{\mathcal{L}}
    
     \newcommand{\HH}{\mathcal{H}}
    \newcommand{\MM}{\mathcal{M}}

 \newcommand{\SC}{\mathcal{S}}

\newcommand{\eps}{\varepsilon}
\newcommand{\g}{\gamma}

\begin{document}

\author{Clementa Alonso-Gonz{\'a}lez}
\address{Universidad de Alicante. Departamento de Matemáticas.
Carretera de San Vicente del Raspeig s/n, E-03690, San Vicente del
Raspeig Alicante (Spain)} \email{clementa.alonso@ua.es}
\author{Fernando Sanz S\'{a}nchez}
\address{Universidad de Valladolid. Departamento de \'{A}lgebra, An\'{a}lisis Matem\'{a}tico, Geometr\'{\i}a y Topolog\'{\i}a. Facultad de Ciencias. Paseo de Belén, 7, E-47011, Valladolid (Spain)} \email{fsanz@uva.es}
\title[Stratification of the dynamics of three-dimensional real flows]
{Stratification of three-dimensional real flows II: A generalization of Poincaré's planar\\ sectorial decomposition}

\maketitle

\begin{abstract}
	Let $\xi$ be an analytic vector field in $\mathbb{R}^3$ with an isolated singularity at the origin and having only hyperbolic singular points after a reduction of singularities $\pi:M\to\R^3$. Assuming certain conditions to be specified throughout the work at hand, we establish a theorem of stratification of the dynamics of $\xi$ that generalizes to dimension three the classical one, coming from Poincar\'{e}, about the decomposition of the dynamics of an analytic planar vector field into {\em parabolic}, {\em elliptic} or {\em hyperbolic} invariant sectors.
\end{abstract}


\footnotetext[1]{The authors were supported by Ministerio de Ciencia e Innovaci\'on (MTM2016-77642-C2-1-P and PID2019-105621GB-I00).}

\section{Introduction}\label{sec:intro}
\begin{figure}[h]
 	\begin{center}
 		\includegraphics[scale=0.5]{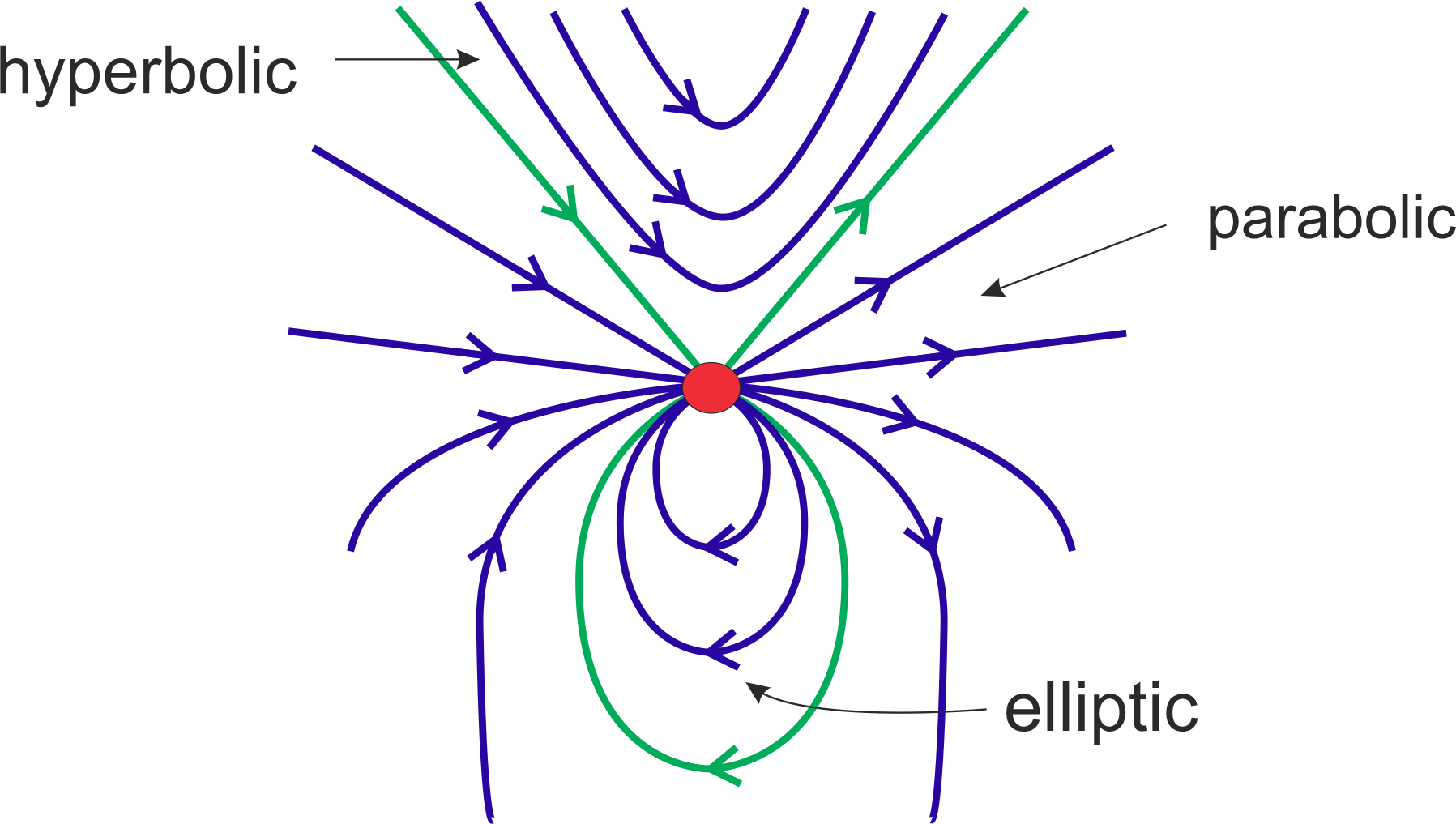}
 		\caption{Sectorial decomposition in the plane.}
 		\label{Fig:PlaneSectors}
 	\end{center}
 \end{figure}

This paper continues our preceding one \cite{Alo-C-C1} on the description of the local dynamics of real analytic three dimensional vector fields. The aim here is to provide a generalization of the classical result about sectorial decomposition of planar vector fields. This result, 
coming from Poincar\'{e} \cite{Poi} and Bendixson \cite{Ben} (see \cite{And} or \cite{Dum} for modern proofs), 
 states that if $\xi$ is a real analytic vector field with an isolated singularity at $0\in\R^2$ and it is not of the {\em center-focus} type, then the dynamics of $\xi$ in a small enough open neighborhood $U$ of $0$ can be decomposed into a finite number of invariant {\em sectors} where all the trajectories of $\xi$ have the same asymptotic behavior (see Figure \ref{Fig:PlaneSectors}): either escaping
 from $U$ in positive and negative time ({\em hyperbolic sector}), either converging to $0$ in both positive and negative infinite time ({\em elliptic sector}), or converging to $0$ in one direction of time and escaping from $U$ in the other one ({\em parabolic sector}). Moreover, the boundary of each sector consists of two trajectories of $\xi$ accumulating at the origin with a defined tangent ({\em characteristic orbits}).
 
The existence of such a sectorial decomposition strongly depends on the geometry of $U$ or, more precisely, in the relative position of the boundary of $U$ with respect to the vector field $\xi$. However, there always exists a base of neighborhoods of $0\in\R^2$ admitting such a decomposition. In any case, the respective numbers of hyperbolic and elliptic sectors do not depend on the neighborhood where a finite sectorial decomposition exists. On the contrary, the number of parabolic sectors is not well determined, even for a fixed neighborhood $U$, since any parabolic sector, being formed by characteristic orbits can be decomposed, in turn, in any number of new parabolic sectors. Moreover, hyperbolic sectors can be ``germified'' in the sense that if we have sectorial decompositions in two respective neighborhoods $U, U'$, there exists a bijection $S\mapsto S'$ between the family of hyperbolic sectors in $U$ and the one in $U'$ such that the germs of $S$ and $S'$ at $0$ coincide. This property does not hold for elliptic sectors (although it is true for the two respective unions of elliptic and parabolic sectors in $U$ and in $U'$).

We want to mention that both, the property that $\xi$ is not of the center-focus type and the configuration of the sectors of any sectorial decomposition depend on a finite jet of the expression of $\xi$ in analytic coordinates (see Dumortier \cite{Dum}). In fact, they can be determined after {\em reduction of singularities} of $\xi$ (in the sense of Seindenberg \cite{Sei}). More precisely, when all resulting singularities are hyperbolic, the structure of the sectorial decomposition can be codified by the configuration of singular points and trajectories of the transformed foliation along the divisor (a directed graph homeomorphic to $\mathbb{S}^1$), as well as by the sign of the transversal eigenvalues at non-corner singular points.

\strut

The problem of sectorial decompositions in higher dimension can be described in terms of local dynamics stratifications.
In general terms, let $\xi$ be a (germ of) real analytic vector field having an isolated singularity at $0\in\R^n$. Let $U$ be an open neighborhood of the origin such that $\xi$ is defined in a neighborhood of $\overline{U}$. A set $A\subset U$ which is saturated by the flow of $\xi|_U$ will be called respectively {\em hyperbolic; elliptic} or {\em  parabolic for $\xi$ (relatively to $U$)} if for any integral curve $\g$ of $\xi|_U$ contained in $A$ one has, respectively (and independently of $\g$):  $\g$ escapes $U$ in both directions;  $\alpha(\g)=\omega(\g)=0$; $\g$ accumulates to $0$ in one direction and escapes from $U$ in the other one. 
We say that $U$ is {\em (analytically) $\xi$-stratifiable} if there exists a finite stratification $\{A_j\}_{j\in J}$ of $U$
into connected analytic submanifolds of $U$ which are either hyperbolic, elliptic or parabolic subsets for $\xi$ relatively to $U$. Recall that a stratification means that $U=\cup_{j\in J}A_j$, that the strata $A_j$ are pairwise disjoint and that, for any stratum $A_j$, the closure of $A_j$ in $U$ is the union of $A_j$ and some other strata $A_{j'}$ with $\dim A_{j'}<\dim A_j$.

Thus, the result about sectorial decomposition of a planar vector field $\xi$ can be reformulated by saying that, when $\xi$ is not of center-focus type, there exists a basis of analytically $\xi$-stratifiable open neighborhoods of the origin.
Coming back to dimension three, our principal contribution  in this paper can be provisionally stated as follows:

\begin{theorem}
	\label{th:beforeblowup-intro}
	Let $\xi$ be an analytic vector field with an isolated singularity at  $0\in\mathbb{R}^3$. Then, under \emph{\textbf{non degeneracy conditions}} over $\xi$, there exists a neighborhood basis $\{U_n\}$ of the origin such that each $U_n$ is analytically $\xi$-stratifiable.
\end{theorem}

The ``non degeneracy conditions'' we impose to $\xi$ are expressed in terms of the resulting foliation after a process of reduction of singularities of the real vector field $\xi$. They have been already considered in the works of Alonso \textit{et al}. \cite{Alo-C-C1,Alo-C-C2,Alo-C-R} as well as in the first part \cite{Alo-San1} in which the paper at hand is based. Let us briefly recall them.

By means of the Panazzolo's work \cite{Pan}, we know that there exists a real analytic proper morphism $\pi:M\to(\R^3,0)$, composition of finitely many {\em weighted blow-ups} with global non-singular centers, where $M$ is a real analytic manifold with boundary and corners, such that, if $\LL_{\xi}$ is the one-dimensional foliation generated by $\xi$ and $\LL=\pi^*\LL_\xi$ is the transformed foliation in $M$, then any singular point $p\in\Sing(\LL)$ is {\em elementary}, i.e., the linear part of a local generator at $p$ is non-nilpotent.
In our case, since $0$ is an isolated singularity of $\xi$, we have, in addition, that the {\em total divisor} $D=\pi^{-1}(0)=\partial M$ is homeomorphic to $\SSS^2$ and $\Sing(\LL)\subset D$ (since in Panazzolo's procedure we just blow-up along centers contained in the singular locus of the intermediate transformed foliations).

In fact, we consider a more general scenario than the one coming as a result of a reduction of singularities of a local vector field: our object of study is a triple tupla $\mathcal{M}=(M,D,\LL)$, called a {\em foliated manifold}, where $M$ is a real analytic manifold with boundary and corners, $D:=\partial M$ is a compact normal crossings divisor and $\LL$ is a one-dimensional oriented singular foliation on $M$ for which $\Sing(\LL)\subset D$. The foliated variety $\mathcal{M}=(M,D,\LL)$ is called:
\begin{itemize}
\item {\em Spherical}, if $D$ is homeomorphic to $\SSS^2$.
\item {\em Non-dicritical}, if $D$ is invariant by $\LL$.
\item {\em Hyperbolic}, if any singularity is hyperbolic (thus $\Sing(\LL)$ is finite since $D$ is compact).
\item {\em Acyclic}, if it is non-dicritical and there are no cycles (closed leaves) nor  polycycles of the restricted foliation $\LL|_D$. This notion can be seen as the generalization of the ``non center-focus type'' restriction to dimension three.
\item {\em Of Morse-Smale type}, if given two different singularities $p,q\in\Sing(\LL)$ which are two-dimensional saddle points for the restriction of $\LL$ to some component of $D$, they cannot be connected by any leaf inside $D$, unless this leaf is contained in the singular locus of $D$ (intersection of two components).
\item  {\em Non saddle-resonant} (\textit{non s-resonant}, for short), if there are no resonant multiple saddle connections.

\end{itemize}

\begin{figure}[h]
	\begin{center}
		\includegraphics[scale=0.6]{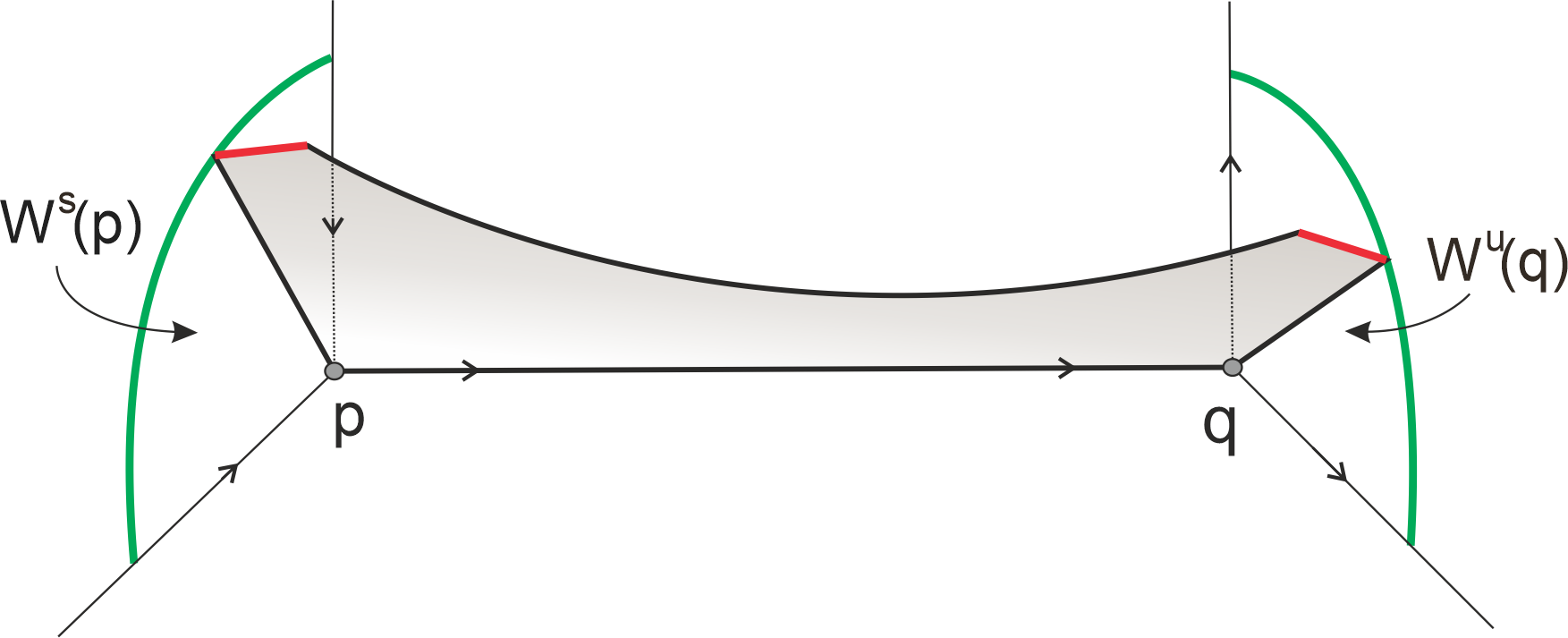}
		\caption{A resonant saddle-connection.}
		\label{Fig:ConexionSillas}
	\end{center}
\end{figure}

The last condition is a little bit more entangled and, just to give an idea of what it means, we can think about a unique connection between two saddle points $p,q$ sharing their one-dimensional invariant manifold that, in turn, are contained in the skeleton. In this case, the non $s$-resonance condition implies that the quotients between the two eigenvalues of the same sign at $p$ and at $q$ are different. 
Dynamically, it means that the situation represented in 
  Figure \ref{Fig:ConexionSillas} is not allowed.

In that picture we can appreciate how the flow saturation of a small curve that accumulates at a midpoint of the two-dimensional invariant manifold at $p$, accumulates also at a midpoint of the corresponding one at $q$. For a complete discussion the reader can consult \cite{Alo-San1} as well as the previous papers \cite{Alo-C-C1, Alo-C-C2, Alo-C-R}. As explained there, in general it concerns non-resonance conditions over the family of quotients of equal-sign-eigenvalues of the linear part of local generators of $\LL$ at some singularities which are three-dimensional saddles. It is satisfied if, for instance, such quotients form a family of algebraically independent real numbers. Our condition is much more weak than that. In any case, it is a generic condition, that is, stable by small perturbations of $\LL$ if the rest of conditions are fulfilled.

Given a foliated manifold $(M,D,\LL)$ and an open neighborhood $U$ of $D$ in $M$, a non-empty subset $A\subset U$ is said to be an {\em elementary dynamical piece in $U$} if $A$ is connected, saturated by $\LL$ inside $U$ (that is, the leaf of $\LL|_U$ at any point of $A$ is entirely contained in $A$) and we can assign to $A$ two sets $\alpha(A),\omega(A)\subset\Sing(\LL)$, each one of them being either empty or unipunctual, such that for any $p\in A$ the leaf $\ell_p\subset A$ of $\LL|_U$ through $p$ satisfies $\alpha(\ell_p)=\alpha(A)$ and $\omega(\ell_p)=\omega(A)$. An elementary dynamical piece will be said  {\em elliptic}, if $\alpha(A)\ne\emptyset\ne\omega(A)$;  {\em hyperbolic}, if $\alpha(A)=\emptyset=\omega(A)$, {\em parabolic attracting}, if $\alpha(A)=\emptyset\ne\omega(A)$ or {\em parabolic repelling}, if $\alpha(A)\ne\emptyset=\omega(A)$. 

Let $k\in\N_{\ge 0}\cup\{\infty,\omega\}$. The open neighborhood $U$ is said to be {\em $\LL$-stratifiable of class $\CC^k$} if there
exists a finite stratification $\mathcal{S}=\{A_j\}_{j\in J}$ of $U$ into locally open connected submanifolds of $U$ of class $\CC^k$ such that each stratum $A_j$ is an elementary dynamical piece of $\LL$ in $U$. Such an $\SC$ is called an {\em $\LL$-stratification} (of class $\CC^k$). As usual, $k=\omega$ means ``real analytic", in which case we talk about an {\em analytic $\LL$-stratification}. 

\strut

The precise statement of the main result that we prove in this article is the following.

\begin{theorem}[Main result]\label{th:afterblowup-intro}
	Let $\mathbb{M}=(M,D,\LL)$ be a foliated manifold which is spherical, non-dicritical, hyperbolic, acyclic, non s-resonant and of Morse-Smale type.  Then there exists a fundamental system $\{U_n\}_{n\in\N}$ of analytic $\LL$-stratifiable open neighborhoods of $D$ in $M$. More precisely, for each $n \in \N$, there exists an analytic stratification $\mathcal{S}_n=\{A^n_j\}_{j\in J_n}$ of $U_n$
	into elementary dynamical pieces in $U_n$, satisfying the following properties:
	\begin{enumerate}[(i)]
		\item $\mathcal{S}_n$ is compatible with the divisor $D$, that is,  for any stratum $A$ of $\mathcal{S}_n$, either $A\cap D=\emptyset$ or $A\subset D$.  Moreover, the restriction $\mathcal{S}_D=\mathcal{S}_n|_D$ is an analytic stratification of $D$ not depending on $n$.
		\item For any $j\in J_n$ and, for any $n\in\N$, we have that $\overline{A^n_j}\cap D$ is a non-empty set (thus a non-empty union of strata of $\mathcal{S}_D$).
		\item The subsets $J_n^H,J_n^E\subset J_n$ of indices corresponding to the hyperbolic and elliptic strata in $\SC_n$, respectively, can be chosen so that they do not depend on $n$ and such that for any $j\in J_n^H\cup J_n^E$, the sets $\alpha(A_j^n),\omega(A_j^n)$ and $\overline{A^n_j}\cap D$ also do no depend on $n$.
		%
		\end{enumerate}

\end{theorem}
\begin{remark}\label{rk:SC-topologic}
	{\em
In the proof of Theorem~\ref{th:afterblowup-intro} we will see that, in fact, we can guarantee item (iii) for the whole subset of indices if we only require the stratification to be of class $\CC^0$. More precisely, up to substituting several parabolic strata in each $\SC_n$ for their union (but preserving the hyperbolic and elliptic strata), we get coarser topological $\LL$-stratifications $\SC'_n=\{A'^n_j\}_{j\in J'}$, where $J'$ is independent of $n$, satisfying (i), (ii) and such that for any $n$ and any $j\in J'$, the character of $A'^n_j$, as well as the sets $\alpha(A'^n_j),\omega(A'^n_j),\overline{A'^n_j}\cap D$ do not depend on $n$. In fact, the difference between $\SC'_n$ and its analytic refinement $\SC_n$ is that in the first one there can be some two-dimensional parabolic strata which are not regular (differentiable) submanifolds. We refine them including the family of curves in the singular locus of these strata in order to get $\SC_n$. The number of components in this singular locus depend on the construction of the elements $U_n$ in the neighborhood basis (see Section~\ref{sec:proof} for further discussion).
}
\end{remark}
Now, it is clear that Theorem~\ref{th:afterblowup-intro} implies Theorem~\ref{th:beforeblowup-intro} if we put as ``non-degeneracy conditions'' the requirement that there exists a reduction of singularities $\pi:M\to\R^3$ of the vector field $\xi$ such that the foliated manifold $(M,\pi^{-1}(0),\pi^*\LL_\xi)$ (already spherical) is hyperbolic, non-dicritical, acyclic, non $s$-resonant and of Morse-Smale type. In fact, the image $\pi(U_n)$ of the set $U_n$ in the statement of Theorem~\ref{th:afterblowup-intro} will provide, in this case,  an open analytically $\xi$-stratifiable neighborhood of the origin with the stratification $\{0\}\cup\{\pi(A_j^n)\,:\,A_j^n\cap D=\emptyset\}$. Notice, moreover, that $\pi(A_j^n)$ is parabolic (attracting or repelling), hyperbolic or elliptic according to the same nature of $A_j^n$ as a stratum in $U_n$. 

\strut

The proof of Theorem~\ref{th:afterblowup-intro} is 
 mainly based on the construction carried out in our first article \cite{Alo-San1} of a neighborhood basis of $D$ formed by {\em fitting domains}. In Section~\ref{sec:background} we discuss this concept and recall the results and ingredients that we need.  Roughly, a fitting domain is a compact neighborhood $F$ of $D$ whose frontier $Fr(F)$ is a piecewise smooth surface  everywhere tangent to the foliation $\LL$ except along some ``entry-exit'' controlled discs $T_p\subset Fr(F)$ associated to points $p$ in some subset $S'$ of $\Sing(\LL)$. Precisely, $S'$ is formed by those singular points for which one of the local invariant manifolds (either stable or unstable), denoted by $W^{tr}(p)$, is transversal to $D$. Such a fitting domain is obtained by a perfect gluing of ``chimney-type'' neighborhoods, tubes and flow boxes associated, respectively, to vertices, edges and faces of a directed planar graph $\Omega=\Omega_\MM$ contained in $D$. This graph schematizes the dynamics of the restriction $\LL|_D$. Its major property, inherited from the acyclicity assumption, is that it has no cycles of edges. Such a property permits to avoid an endless recurrence in the gluing process, which is crucial in order to get fitting domains.

 In Section~\ref{sec:proof} we prove Theorem~\ref{th:afterblowup-intro} by showing that the interior $U=\text{int}(F)$ of a fitting domain $F$ has a finite analytic stratification into elementary dynamical pieces. Vertices, edges and faces of the graph $\Omega$ already give a stratification $\mathcal{S}_D$ of $D$  independent of $U$. The structure of the strata outside the divisor is related to the family $\{W^{tr}(p)\}_{p\in S'}$ of transversal manifolds mentioned above. 
 Since $W^{tr}(p)$ is composed of leaves of $\LL$ accumulating to $p$ in only one sense, the $\LL$-saturation of $W^{tr}(p)\cap\text{Fr}(F)$ (a subset of $T_p$)  intersects $U$ along parabolic strata of dimension at most $\dim(W^{tr}(p))$. These are the local strata attached to the singular point $p\in S'$ and its number (unless we only require a topological stratification) may depend on the number of singular points in the boundary of the disc $T_p$. 
 
 The main global strata in the complement of the local strata are generated at those $p\in S'$ for which $\dim(W^{tr}(p))=2$, the so called {\em transversal saddles}. 
  More precisely, attached to any transversal saddle $p$, we have two types of two-dimensional strata. In one hand, the saturation of $W^{tr}(p)$ by the flow of $\LL$ ``separates'' the dynamics and must be considered as part of the union of two-dimensional strata. It is called the {\em fixed separating surface} generated at $p$. Its germ along $D$ does not depend on the fitting domain $U$. On the other hand, if $T_p\subset\text{Fr}(F)$ is the transversal disc mentioned above, there are several intervals inside the boundary $\partial T_p$ whose saturation ``enters'' $U$ (in one of the two senses of the flow), separating also the dynamics. They must also be included in the union of two-dimensional strata. We call them {\em mobile separating surfaces}, since they depend on the chosen fitting domain $F$. 
  
  In order to have a stratification with a finite number of strata, we use a special extra property of the fitting domains from the construction in \cite{Alo-San1}. Notably, that two elements in the family of fixed and mobile separating surfaces do not intersect, as long as the fitting domain is small enough. Morally, the separating surfaces have ``non-oscillating properties''. Moreover, the closure of a fixed or mobile separating surface intersects the divisor only along a path of edges of the graph $\Omega$.
  This will permit to prove items (i) and (ii) of Theorem~\ref{th:afterblowup-intro}. Item (iii) will be a consequence of the construction and the fact that the number of fixed and mobile separating surfaces is independent of the fitting domain.

In Section~\ref{sec:miscellaneous} we consider miscellaneous examples and discuss some questions related to the result of stratification. In particular, we deal with the problem of whether the union of hyperbolic strata is ``germifiable'', as it was the case for the sectorial decomposition of planar vector fields. We provide a concrete example showing that this property is not longer true in dimension three.

\section{Background on spherical foliated manifolds}\label{sec:background}
In this section we summarize the content of the paper \cite{Alo-San1} that we need in order to prove Theorem~\ref{th:afterblowup-intro}. 
\subsection{Elements of a spherical hyperbolic foliated variety}
We fix a  foliated manifold $\MM=(M,D,\LL)$, where $M$ is a three-dimensional analytic manifold with boundary and corners, $D=\partial M$ and $\LL$ is an oriented one-dimensional analytic foliation in $M$. Along the paper, we use basic and standard notations concerning leaves and saturation of a set with respect to $\LL$. Namely, if $a\in M$, we denote, respectively, by $\ell_a,\ell^+_a$ and $\ell^-_a$, the leaf, positive leaf and negative leaf of $\LL$ through $a$. If $A\subset M$ and $a\in A$, $\ell^{A}_a$ denotes the connected component of $\ell_a\cap A$ containing $a$ and it is called the {\em $A$-leaf} through the point $a$. Also, if $B\subset A\subset M$, the {\em saturation of $B$ in $A$} is the set $\Sat_A(B)=\bigcup_{b\in B}\ell^A_b$. Replacing $\ell_a$ by $\ell^+_a$ or $\ell^-_a$ everywhere, we define analogously the positive saturation $\Sat^+_A(B)$ or negative saturation $\Sat^-_A(B)$. When $U\subset M$ is an open set and $\ell$ is an $U$-leaf, or a positive $U$-leaf, we denote by $\omega_U(\ell)$ its $\omega$-limit set (inside $U$) as a leaf of the restricted foliation $\LL|_U$. For instance, if $\ell$ is contained in some leaf $\tilde{\ell}$ that cuts the frontier of $U$ then we have $\omega_U(\ell)=\emptyset$ although it is possible that $\omega_M(\tilde{\ell})\ne\emptyset$. We delete the index ``$U$'' when $U=M$.

Using the definitions in the introduction, we assume from now on that $\MM=(M,D,\LL)$ is spherical, non-dicritical, hyperbolic and acyclic; that is, 
$D$ is homeomorphic to the sphere $\SSS^2$, $D$ is invariant for $\LL$ (i.e., $\Sat_M(D)=D$), every singular point $p\in\Sing(\LL)$ belongs to $D$ and it is hyperbolic (eigenvalues with non-zero real part) and, moreover, the restriction $\LL|_D$ has no cycles nor polycycles. 

The boundary $D$ is called the {\em divisor}. 
A {\em component of $D$} is the closure of a connected component of the set $\Reg D\subset D$ of points where $D$ is smooth. The function $e:D\to\{1,2,3\}$ which assigns to any $a\in D$ the number $e(a)$ of components passing through $a$ provides a stratification of $D$, called the {\em standard stratification} and denoted by $Stan(D)$. The invariance condition implies that any stratum of $Stan(D)$ is an analytic smooth submanifold of $M$, tangent everywhere to $\LL$. In particular, the $0$-dimensional strata (corners) are singular points of $\LL$ and the $1$-dimensional strata are non-singular leaves of $\LL$. The {\em skeleton} of $D$, denoted by $Sk(D)$, is the union of strata in $Stan(D)$ of dimension less or equal to 1; i.e., the skeleton is the complementary of $\Reg(D)$ in $D$, or the set of points where at least two components of $D$ meet.

Given $p\in\Sing(\LL)$, denote by $W^s(p)$ and $W^u(p)$ respectively the local stable and unstable manifold of a local generator of $\LL$ at $p$. By the condition that the divisor is invariant, we have that at least one of them is contained in $D$. The point $p$ is called a {\em tangential saddle} if $W^s(p)\cup W^u(p)\subset D$ (notice that necessarily in this case $p$ is a three-dimensional saddle, i.e., both $W^s(p)$ and $W^u(p)$ are of positive dimension). Denote by $S_{tg}\subset\Sing(\LL)$ the set of tangential saddles and by $S'=\Sing(\LL)\setminus S_{tg}$ (as in the introduction). If $p\in S'$, denote by $W^{tr}(p)$ the element of $\{W^s(p),W^u(p)\}$ which is not contained in $D$. In this case, if $\dim W^{tr}(p)=1$ or $\dim W^{tr}(p)=3$ then $p$ is called a {\em $D$-node} (the name is justified because, in both cases, the restriction of $\LL$ to any component of $D$ has a node singularity at $p$); if $\dim W^{tr}(p)=2$ then $p$ is called a {\em transversal saddle}. Denote by $N$ the set of $D$-nodes and by $S_{tr}$ the set of transversal saddles, so that we have the partition $\Sing(\LL)=S_{tg}\cup N\cup S_{tr}$. We would need to consider separately all the possibilities for points in $S'=N\cup S_{tr}$, so that we introduce the notation
$$
N=N^s_1\cup N^u_1\cup N^s_3\cup N^u_3,\;\;S_{tr}=S_{tr}^s\cup S_{tr}^u,
$$ 
where, for $a\in\{s,u\}$ and $i\in\{1,3\}$, $p\in N^a_i$ if $W^{tr}(p)=W^a(p)$ and $\dim W^{tr}(p)=i$ and, for $a\in\{s,u\}$, $p\in S_{tr}^a$ if $W^{tr}(p)=W^a(p)$ and $\dim W^{tr}(p)=2$. We put also $N_i:=N_i^s\cup N_i^u$ for $i=1,3$.

\subsection{Limit sets of leaves}
Let us state a first dynamical consequence for the foliated manifold $\MM=(M,D,\LL)$. 

\begin{proposition}\label{pro:omega-limits}
Suppose that $\MM$ is spherical, non-dicritical, hyperbolic and acyclic. Then there exists an open neighborhood $U$ of $D$ in $M$ with the following property: if $\ell^+$ is a positive $U$-leaf,  then either $\omega_U(\ell^+)=\emptyset$ or $\omega_U(\ell^+)$ consists of a unique point of $\Sing(\LL)$ (in this last case $\ell^+$ is also a positive $M$-leaf and $\omega_U(\ell^+)=\omega(\ell^+)$). The same conclusion holds for the $\alpha$-limit set of any negative $U$-leaf.
\end{proposition}
\begin{proof}
We prove first that there exists an open neighborhood $U$ such that no point of $U\setminus D$ can be an $\omega$-limit point of any positive $U$-leaf. Suppose, on the contrary, that there exists a fundamental system $\{U_n\}_n$ of open neighborhoods of $D$, where each $U_n$ has a positive $U_n$-leaf $\tau_n\subset U_n$ with the property that $\omega_{U_n}(\tau_n)\cap(U_n\setminus D)\ne\emptyset$. Notice that $\tau_n$ is also a positive leaf of $M$ and $\omega_{U_n}(\tau_n)=\omega(\tau_n)\cap U_n$ since, otherwise, $\tau_n$ would cross the frontier $Fr(U_n)$ and hence $\omega_{U_n}(\tau_n)$ would be empty. In particular, $\omega(\tau_n)$ is closed in $M$ and invariant by $\LL$ (that is, if $x\in\omega(\tau_n)$ then $\ell_x\subset\omega(\tau_n)$). Let $Z$ be the set consisting of those points $a\in M$ for which there is a sequence $\{b_{k}\}_{k}$ with $b_{k}\in\omega(\tau_{n_k})\cap(U_{n_k}\setminus D)$ for any $k$ and $\lim_{k\to\infty} n_k=\infty$ such that $a=\lim_{k\to\infty}b_{n_k}$.
Notice that $Z$ is a non-empty closed set contained in $D$ and invariant by $\LL$, since any $\omega(\tau_n)$ is invariant. 
We conclude, using Poincaré-Bendixson's Theorem on $D\simeq\SSS^2$, together with the acyclicity condition, that $Z$ must contain at least one singular point $p$. 

Notice that $p$ cannot belong to $N_3^u\cup N^s_3$. Otherwise, there would be a neighborhood $B$ of $p$ such that no point in $B\setminus\{p\}$ can be the $\omega$-limit point of a positive leaf in $M$ (if $x\in B\setminus\{p\}$, the positive leaf $\ell^+_x$ at $x$ either escapes from $B$, when $p\in N^u_3$, or accumulates only in $p$, when $p\in N^s_3$). 

Hence $p$ is a three-dimensional saddle point. We claim that $Z$ intersects both invariant manifolds $W^u(p)$ and $W^s(p)$ at points different from $p$. To show this, put $p=\lim_{k\to\infty} b_{k}$ where $b_{k}\in\omega(\tau_{n_k})\cap(U_{n_k}\setminus D)$ and let $\ell_{b_k}$ be the leaf through $b_{k}$ for any $k$. Notice that $\ell_{b_k}\subset\omega(\tau_{n_k})$ and $\ell_{b_k}$ does not intersect neither $W^u(p)$ nor $W^s(p)$, by a similar argument as the one used to prove that $p\not\in N_3^u\cup N^s_3$. Using Hartman-Grobman's Theorem, or the $\lambda$-lemma (see \cite{Pal-M} for instance), there exists a compact neighborhood $V$ of $p$ such that the sequence $\{\ell_{b_k}\cap V\}_k$ of compact sets accumulates to a set containing points of $W^u(p)\setminus\{p\}$ and of $W^s(p)\setminus\{p\}$. 

Taking into account this claim, we must have $p\not\in N_1^u\cup N_1^s$ and $W^u(p)\subset D$. Choose some point $a_1\in Z\cap W^u(p)\setminus\{p\}$ and let $\sigma_1$ be the leaf through $a_1$. We have that $\sigma_1$, together with the point $p_1=\omega(\sigma_1)$, are contained in $Z$. Considering the same reasoning with $p_1$ and repeating the argument, we create a sequence of non-singular leaves $\sigma_1,\sigma_2,...$ in $D$ with the property that $\alpha(\sigma_{j+1})=\omega(\sigma_j)$ for any $j\ge 1$. Since $\Sing(\LL)$ is finite, such a sequence must contain a polycycle, which is prohibited by the acyclicity hypothesis.

We conclude that there exists an open neighborhood $U$ of $D$ such that, given a positive $U$-leaf $\ell^+$, either $\omega_U(\ell^+)=\emptyset$ or $\emptyset\ne\omega_U(\ell^+)=\omega(\ell^+)\subset D$. To finish, we must show that, in this last case, $\omega(\ell^+)$ reduces to a singular point. By Poincaré-Bendixson's Theorem and the acyclicity condition, this is the case if $\ell^+\subset D$. Suppose that $\ell^+\subset U\setminus D$ and that $\omega(\ell^+)\ne\emptyset$. Again, since there are no cycles, we have that $\omega(\ell^+)$ contains a singular point $p$. If $\ell^+\cap W^s (p)\ne\emptyset$, then we have $\omega(\ell^+)=\{p\}$ as wanted. On the other hand, if $\ell^+\subset W^u(p)$ in a given neighborhood of $p$ then $p\not\in\omega(\ell^+)$. We may assume hence that $\ell^+$ cuts the complement of $W^u(p)\cup W^s(p)$ in any neighborhood of $p$ arbitrarily small. Using similar arguments as above, we construct a non-trivial polycycle contained in $\omega(\ell^+)$, again impossible due to the acyclycity condition.
\end{proof}

\subsection{The associated graph}
We associate a planar oriented graph $\Omega=\Omega_\MM$ with $\MM=(M,D,\LL)$ as follows:  
\begin{itemize}
	\item The set of vertices is $V(\Omega):=\Sing(\LL)$.
	\item The set of edges $E(\Omega)$ is the family  of non-singular leaves $\sigma$ of $\LL$ inside $D$ satisfying one of the following properties: either $\sigma$ is a one-dimensional stratum of $Stan(D)$, or there is a component $D_i$ of $D$ and some point $p\in D_i\cap\Sing(\LL)$ such that $p$ is a two-dimensional saddle point of the restriction $\LL|_{D_i}$ and $\sigma$ is the leaf containing one of the two local invariant manifolds, stable or unstable, of $\LL|_{D_i}$ at $p$. The orientation of an edge is that induced by the orientation of $\LL$.
	\item A vertex $p$ is adjacent to an edge $\sigma$ iff $p=\alpha(\sigma)$ or $p=\omega(\sigma)$ (notice that the $\alpha$ or $\omega$-limit set of a leaf in $D$ is a singleton by Poincaré-Bendixson's Theorem and the acyclicity condition).
\end{itemize}

The support of $\Omega$ is the set $|\Omega|\subset D$ formed by the union of edges and vertices.

Let $p,q\in V(\Omega)$. Recall that a {\em path of edges} of $\Omega$ {\em from $p$ to $q$} is a sequence $\g=(\sigma_0,\sigma_1,...,\sigma_{n})$ where each $\sigma_i$ is an edge of $\Omega$ going from vertices $p_i$ to $p_{i+1}$ such that $p_0=p$ and $p_{n+1}=q$. The support of $\g$, denoted by $|\g|$ is the union of the edges $\sigma_i$ and their adjacent vertices. The acyclicity condition implies also that there are no non-trivial paths from a vertex to itself (no cycles in the graph).

\subsection{Additional hypothesis}\label{sec:additional}
As mentioned in the introduction, we impose two more conditions to $\MM=(M,D,\LL)$.

 {\bf Morse-Smale type.-}  A {\em saddle connection} is an edge $\sigma\in E(\Omega)$ for which there exists a component $D_i$ of $D$ containing $\sigma$ and such that both extremities $p=\alpha(\sigma)$ and $q=\omega(\sigma)$ are two-dimensional saddle points of the restriction $\LL|_{D_i}$. A {\em multiple saddle connection} is a path of edges consisting of saddle connections. We say that $\MM=(M,D,\LL)$ is of {\em Morse-Smale type} if any saddle connection is contained in the skeleton.

{\bf No saddle-resonances.-} This condition needs the introduction of some terminology. Let us recall here briefly this notion in a mere combinatoric way, but see \cite{Alo-C-C1, Alo-San1} for the dynamical interpretation of it.   

Let $\sigma$ be an edge contained in the skeleton and denote by $\CC_\sigma$ the set of components of the divisor where $\sigma$ is contained (notice that $\CC_\sigma$ consists of two elements). A {\em weight assignment} of $\sigma$ is a map
$$
\rho_\sigma:\CC_\sigma\to\R_{>0}
$$
such that, if $\CC_{\sigma}=\{D_1,D_2\}$ then $\rho_{\sigma}(D_1)\rho_{\sigma}(D_2)=1$. Notice that $\rho_\sigma$ is completely determined by only one of the values $\rho_\sigma(D_1)$ or $\rho_\sigma(D_2)$.

If $p\in\Sing(\LL)$ is a three-dimensional saddle point, we denote by $W^i_p$, for $i\in\{1,2\}$, the local $i$-dimensional invariant manifold of $\LL$ at $p$. Denote  $\alpha(p),\lambda(p),\lambda'(p)$ the three eigenvalues of the linear part of a local generator of $\LL$ at $p$, where $\alpha(p)$ is the one associated with $W^1_p$ and $\{\lambda(p),\lambda'(p)\}$ the associated with $W^2_p$, that is, $\lambda(p),\lambda'(p)$ have real parts of the same sign (notice that these eigenvalues are well defined up to a common real positive multiple). 
Suppose that $W^1_p$ is contained in the skeleton of $D$ and let $\sigma$ be one edge intersecting $W^1_p$. Write also   $\CC_\sigma=\{D_1,D_2\}$ and denote by $\tau_i$ the edge of $\Omega$ containing $D_i\cap W^2_p$, for $i\in\{1,2\}$. Notice that the tangent of $\tau_i$ at $p$ is an eigendirection of the linear part of $\LL$, with associated eigenvalue $\lambda_i$ (necessarily a real number) so that $\{\lambda_1,\lambda_2\}=\{\lambda(p),\lambda'(p)\}$. We define in this case a weight assignment $w_{\sigma,p}:\CC_\sigma\to\R_{>0}$ of $\sigma$, called {\em $p$-structural}, by 
$$
w_{\sigma,p}(D_i):=\frac{\lambda_j}{\lambda_i},\;\;\;\mbox{ where }\{i,j\}=\{1,2\}.
$$
In case $W^1_p$ is contained in the skeleton, and using the same notations, given a weight assignment $\rho_\sigma$ of $\sigma$ different from the $p$-structural one, we can ``transport'' it, either to $\tau_1$ or to $\tau_2$, as follows. If $i\in\{1,2\}$ is the (unique) index such that $\rho_\sigma(D_i)>w_{\sigma,p}(D_i)$ then we say that there exists the {\em transition weight of $\rho_\sigma$ from $\sigma$ to $\tau_i$}, to be defined as the weight assignment $\rho_{\tau_i}=T_{\sigma,\tau_i}(\rho_\sigma)$ of $\tau_i$ given by
\begin{equation}\label{eq:transition}
	\rho_{\tau_i}(D_i)=\frac{\lambda_j-\lambda_i \rho_{\sigma}(D_i) }{\alpha(p)}.
\end{equation}
Reciprocally, given a weight assignment $\rho_{\tau_i}$ of $\tau_i$, where $i \in \{1,2\}$, the {\em transition of $\rho_{\tau_i}$ from $\tau_i$ to $\sigma$} is the weight assignment $\rho_\sigma=T_{\tau_i,\sigma}(\rho_{\tau_i})$ of $\sigma$ defined by setting the value $\rho_{\sigma}(D_i)$ that satisfies (\ref{eq:transition}). It turns out  that this assignment $\rho_\sigma$ is different from the $p$-structural weight assignment, that there exists the transition of $\rho_\sigma$ from $\sigma$ to $\tau_i$ and that $T_{\sigma,\tau_i}(\rho_{\sigma})=\rho_{\tau_i}$; that is, the transition $T_{\tau_i,\sigma}$ is the ``inverse'' of the transition $T_{\sigma,\tau_i}$.

Let $p,q\in\Sing(\LL)$ be two three-dimensional saddles with $p\ne q$ and let $\g=(\sigma_0,\sigma_1,...,\sigma_n)$ be a multiple saddle connection from $p$ to $q$. Denote  $(p_0,p_1,...,p_{n+1})$ the corresponding sequence of consecutive vertices being $p_0=p$ and  $p_{n+1}=q$. We say that $\g$ is {\em saddle-resonant} (or just {\em s-resonant}) if the following properties hold:
\begin{itemize}
	\item Each edge $\sigma_i$ is contained in the skeleton.
	\item The first edge $\sigma_0$ intersects $W^1_p$ and the last edge $\sigma_{n}$ intersects $W^1_q$.
	\item Denote by $\rho_0=w_{p,\sigma_0}$ the $p_0$-structural weight assignment of $\sigma_0$. Then we have $\rho_0\ne w_{p_1,\sigma_0}$ and there exists the transition of $\rho_0$ from $\sigma_0$ to $\sigma_1$.
	\item Defining recursively
	$
	\rho_i:=
			T_{\sigma_{i-1},\sigma_{i}}(\rho_{i-1}),
	$
	 for $i=1,...,n-1$, we have that $\rho_{i}\ne w_{p_{i+1},\sigma_{i+1}}$ and there exists the transition of $\rho_{i}$ from $\sigma_{i}$ to $\sigma_{i+1}$.
	\item It holds $\rho_{n}:=T_{\sigma_{n-1},\sigma_{n}}(\rho_{n-1})=w_{q,\sigma_n}$.
\end{itemize}
\begin{remark}
\emph{Our hypothesis of non s-resonance is just the requirement that there are no s-resonant multiple saddle connections. Taking into account the polynomial expression in (\ref{eq:transition}) for the transition of weight assignments, we can guarantee this hypothesis if, for instance, the family of quotients $\{Re(\lambda(p))/\alpha(p),Re(\lambda'(p))/\alpha(p)\}$ is algebraically transcendent over $\Q$, where $p$ runs over the set of three-dimensional saddle singular points of $\LL$.} 
\end{remark}

From now on, together with the rest of already stated conditions, we assume that $\MM=(M,D,\LL)$ is of Morse-Smale type and has no s-resonant multiple saddle connections.
One first important dynamical consequence of these hypothesis is the following result (see \cite[Thm. 15 and Cor. 17]{Alo-San1}), that will be crucial in what follows. Denote by $\Omega=\Omega_\MM$ the associated graph of $\MM$.

\begin{proposition}\label{pro:saturation-marks}
 Let $p\in S_{tr}$ be a transversal saddle and denote by $\nu_p^+,\nu_p^-$ the two connected components in which $W^2_p$ divides a (sufficiently small) neighborhood of $p$. Then:
 \begin{enumerate}[(i)]
 	\item There are two paths of edges $\Delta_p^1,\Delta_p^2$ in $\Omega$, each of them with extremities at $p$ and a $D$-node point, satisfying that $$\overline{\Sat_U(W^2_p)}\cap D=|\Delta_p^1|\cup|\Delta_p^2|$$ for any sufficiently small neighborhood $U$ of $D$ in $M$.
 	\item For each $\epsilon\in\{+,-\}$, there is a path $\Theta(\nu_p^\epsilon)$ in $\Omega$, with extremities at $p$ and a $D$-node point, satisfying the following property: for any sufficiently small neighborhood $U$ of $D$ in $M$, if $J$ is a connected parametrized curve contained in $\nu_p^\epsilon\cap U$ such that $\overline{J}\cap  W^2_p=\{a\}$ with $a\in W^2_p\setminus D$, then we have $\overline{\Sat_U(J)}\cap D=|\Theta(\nu_p^\epsilon)|$.
 \end{enumerate}
\end{proposition} 
We notice that the paths $\Delta_p^1,\Delta_p^2$ in item (i) start (resp. end) at $p$ when $W^2_p$ is the unstable (resp. stable) manifold. Also, the paths $\Theta(\nu^+_p),\Theta(\nu^-_p)$ in item (ii) start (resp. end) at $p$ when $W^2_p$ is the stable (resp. unstable) manifold.

\subsection{Fitting domains}\label{sec:fitting}
Our proof of Theorem~\ref{th:afterblowup-intro} is based on the existence of the so called {\em fitting domains}, proved in \cite{Alo-San1}. These are compact neighborhoods of the divisor $D$ whose boundary has very specific properties of transversality/tangency with respect to the foliation $\LL$.  We discuss here these properties since they are needed in the sequel.

Let $B$ be a compact neighborhood of $D$ and assume that $B$ is a semi-analytic subset of $M$. Notice that $Fr(B)\cap\Sing(\LL)=\emptyset$ since we have assumed that the singularities of $\LL$ are contained in $D$. Given $a\in Fr(B)$, by means of the semi-analyticity of $B$ along with the analyticity of $\LL$, we have that each of the two half (open) leaves $\ell^-_a\setminus\{a\},\ell^+_a\setminus\{a\}$ is  contained, either in $int(B)$, or in $ext(B)=M\setminus B$, or in $Fr(B)$, locally at the point $a$. 
Let $\sigma:\{i,e,t\}\to\{int(B),ext(B),Fr(B)\}$ be defined by $\sigma(i)=int(B)$, $\sigma(e)=ext(B)$, $\sigma(t)=Fr(B)$ and say that the point $a\in Fr(B)$ is of {\em type u-v relatively to $B$}, where $u,v\in\{i,e,t\}$, if $\ell^-_a\setminus\{a\}\subset\sigma(u)$ and $\ell^+_a\setminus\{a\}\subset\sigma(v)$. Denote by $Fr(B)^\smallsmile$ the topological interior inside $Fr(B)$ of the set of points of type t-t relatively to $B$. In other words, $Fr(B)^\smallsmile$ is the set of points of $Fr(B)$ where $Fr(B)$ is locally invariant for $\LL$.

In the following definition we use the notations introduced in Section~\ref{sec:intro}. Notably, for each $p\in S'=N\cup S_{tr}$ (a $D$-node or a transversal saddle point) we denote by $W^{tr}(p)$ the invariant manifold at $p$ not contained in $D$. To be precise, only the germ of $W^{tr}(p)$ is well defined, but we fix a representative of it in some given neighborhood of $p$ such that its closure $\overline{W^{tr}(p)}$ is an analytic manifold with boundary and corners homeomorphic to $[0,1]^{\dim W^{tr}(p)}$ and semianalytic in $M$  (see \cite{Car-S},  for instance, for a proof of the analyticity of $W^{tr}(p)$). When $p\in S_{tr}$, will also assume that $\partial W^{tr}(p)\cap D$ is connected and that $\partial W^{tr}(p)\setminus D$ is a smooth curve everywhere transversal to the foliation (where $\partial W^{tr}(p)$ denotes the boundary of $W^{tr}(p)$ as a topological manifold with boundary).

In what follows, for simplicity, 
an {\em (open or closed) interval} (resp. a {\em disc}) is any subset $H\subset M$ homeomorphic to the interval $I=(0,1)$ or $I=[0,1]$ of $\R$ (resp. to $I=\mathbb{D}$, the unit disc of $\R^2$). In both cases, we put $int(H)$, $\partial H$ to denote, accordingly, the subset sent by such a homeomorphism into $int(I)$, $\partial I$. 

\begin{figure}[h]
	\begin{center}
		\includegraphics[scale=0.5]{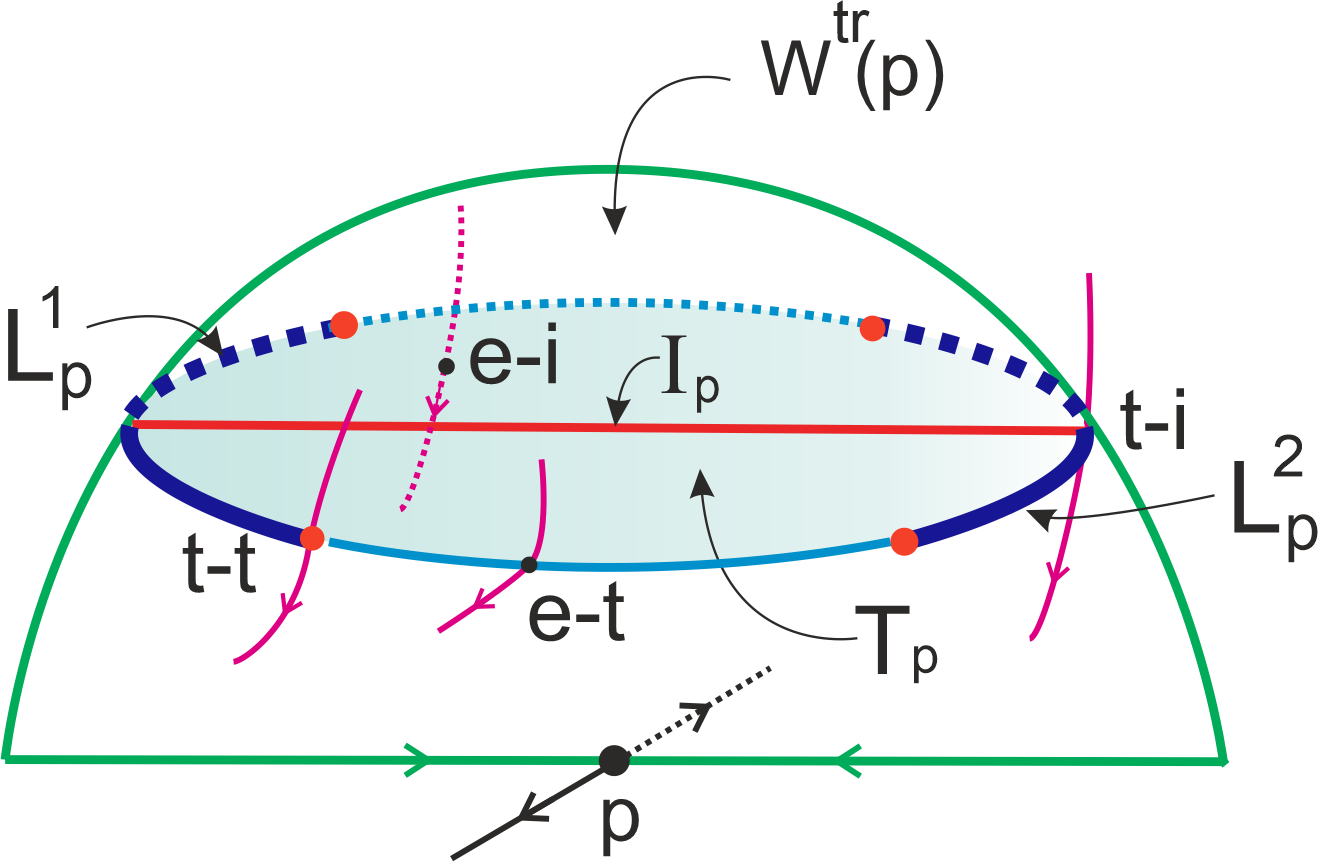}
		\caption{The disc $T_p$ and associated elements if  $\dim W^{tr}(p)=2$.}
		\label{Fig:Disc T_p}
	\end{center}
\end{figure}

\begin{definition}\label{def:fitting}
  A {\em fitting domain} of the foliated manifold $\MM=(M,D,\LL)$ is a compact semianalytic neighborhood $F$ of $D$ in $M$ whose frontier $Fr(F)$ is a topological, piecewise smooth surface and such that there is a family $\{T_p\}_{p\in S'}$ of mutually disjoint non-empty semi-analytic discs satisfying  the following:
\begin{enumerate}[(i)]
	\item  The frontier $Fr(F)$ is the disjoint union
	of $Fr(F)^\smallsmile$ and $\,\bigcup_{p\in S'}T_p$.
	\item Each disc $T_p$ contains $Fr(F)\cap \overline{W^{tr}(p)}$ and  $\overline{W^{tr}(p)}\cap T_p$ coincides with: the disc $T_p$, when $\dim(W^{tr}(p))=3$; a singleton in $int(T_p)$, when $\dim(W^{tr}(p))=1$; a closed non-empty interval $I_p$ such that $int(I_p)\subset int(T_p)$ having extremities in $\partial T_p$ , when $p\in S_{tr}$. 
	\item The type, relatively to $F$, of the points in $T_p$ is as follows (in case $W^{tr}(p)$ is the stable manifold, otherwise, put v-u instead of  u-v):
	\begin{itemize}
		\item Points in $int(T_p)$ are of type e-i.
		\item If $\dim W^{tr}(p)=3$, points in $\partial T_p$ are of type t-i.
		\item If $\dim W^{tr}(p)=1$, points in $\partial T_p$ are of type e-t.		
		\item  If $\dim W^{tr}(p)=2$, there are exactly four points in $\partial T_p$ of type t-t, none of them in $I_p$. Among the four open intervals in which these points divide $\partial T_p$, there are two of them, say $L_p^1$ and $L_p^2$, intersecting $I_p$ and formed by points of type t-i, while the other two do not intersect $I_p$ and are formed by points of type e-t (see Figure \ref{Fig:Disc T_p}).
	\end{itemize}
   \item The elements in the family 
   $
   	\{\Sat_F(\overline{W^{tr}(p)}\setminus D)\}_{p\in S_{tr}}\,\cup\,\{\Sat_F(\overline{L^i_p}\setminus I_p))\}_{p\in S_{tr},i=1,2}
   $ 
   are mutually disjoint subsets of $F$. 
\end{enumerate}
\end{definition}

The main result in \cite{Alo-San1} is then stated as follows.
\begin{theorem}\label{th:fitting}
Assume that $\MM=(M,D,\LL)$ is spherical , non-dicritical, hyperbolic, acyclic of Morse-Smale type and non s-resonant. Then there exists a neighborhood basis $\{F_n\}_{n}$ of $D$ in $M$ such that each $F_n$ is a fitting domain.
\end{theorem}

\section{Proof of the main result}\label{sec:proof}

This section is devoted to prove Theorem~\ref{th:afterblowup-intro}. 
From now on, we consider a foliated manifold $\MM=(M,D,\LL)$ which is spherical, non-dicritical, hyperbolic, acyclic and of Morse-Smale type  without s-resonances. Fix also transversal invariant manifolds $W^{tr}(p)$ for any $p \in S_{tr}$ with the properties discussed in Paragraph~\ref{sec:fitting}. In light of Theorem~\ref{th:fitting}, we can consider a neighborhood basis $\{F_n\}_n$ of $D$ given by fitting domains. Let $U_n:=int(F_n)$ for each $n$, where $int$ means topological interior in $M$. Let us prove that the given neighborhood basis $\{U_n\}_n$ of $D$ satisfies all the requirements in the statement of Theorem~\ref{th:afterblowup-intro}.  

For the sake of simplicity, we drop the subscript $n$ and we simply write $F:=F_n$  and $U:=U_n=int(F_n)$ to denote the fitting domain and its interior, respectively.

\subsection{Fixed and mobile separating surfaces}

 Let us define first some objects, the separating surfaces, whose properties are crucial to determine the three-dimen\-sional strata and, at the same time, contain the more significant strata of dimension two. They are obtained as the saturation by $\LL$ of some of the elements that appear in the definition of fitting domain (cf. Definition~\ref{def:fitting}). We keep the same notations that appear in that definition throughout this paragraph. If $p\in S_{tr}$ and $T_p\subset Fr(F)\setminus Fr(F)^\smallsmile$ is the transversal disc associated to $p$, recall that $L_p^1, L^2_p$ are the two open intervals consisting of points in $\partial T_p$ of type i-t (when $p\in S_{tr}^u$) or of type t-i (when $p\in S_{tr}^s$), relatively to $F$ (see Figure \ref{Fig:Disc T_p}). Each such interval $L^i_p$ intersects $W^{tr}(p)$ at exactly one point, denoted by $a^i_p$, which belongs to the interval $I_p$. For $i=1,2$, let $L^{i,+}_p, L^{i,-}_p$ be the two connected components of $L^i_p\setminus\{a^i_p\}$, again open intervals.

\begin{definition}
Given $p\in S_{tr}$, the set $\HH_p:=\Sat_{U}(W^{tr}(p)\cap U\setminus D)$ is called the {\em fixed separating surface generated at $p$} (on the fitting domain $F$). Also, each of the sets $\GG_p^{i,\epsilon}:=U\cap\Sat_F(L^{i,\epsilon}_p)$, where $i\in\{1,2\}$ and $\epsilon\in\{+,-\}$ is called a {\em mobile separating surface generated at $p$} (on the fitting domain $F$).
\end{definition}
\begin{remark}\label{rk:disjoint-separating}
{\em Notice that each separating surface (fixed or mobile) is contained in one of the elements of the family described in item (iv) of   Definition~\ref{def:fitting}. In particular, any pair of different separating surfaces has an empty intersection.
	
}
\end{remark}

\subsection{Fixed separating surfaces properties}\label{subsec:prop-fixed-separating}
Let $\HH:=\HH_p$ be the fixed separating surface generated at some $p\in S_{tr}$. By simplicity, throughout this subsection we put $W:=W^{tr}(p)$. Hence, by definition, $\HH=\Sat_U(W\cap U\setminus D)$. 
We gather the main properties of $\HH$ by assuming that  $W$ is the stable manifold at $p$. The corresponding ones in case it is the unstable manifold hold analogously.

	\begin{figure}[h]
	\begin{center}
		\includegraphics[scale=0.6]{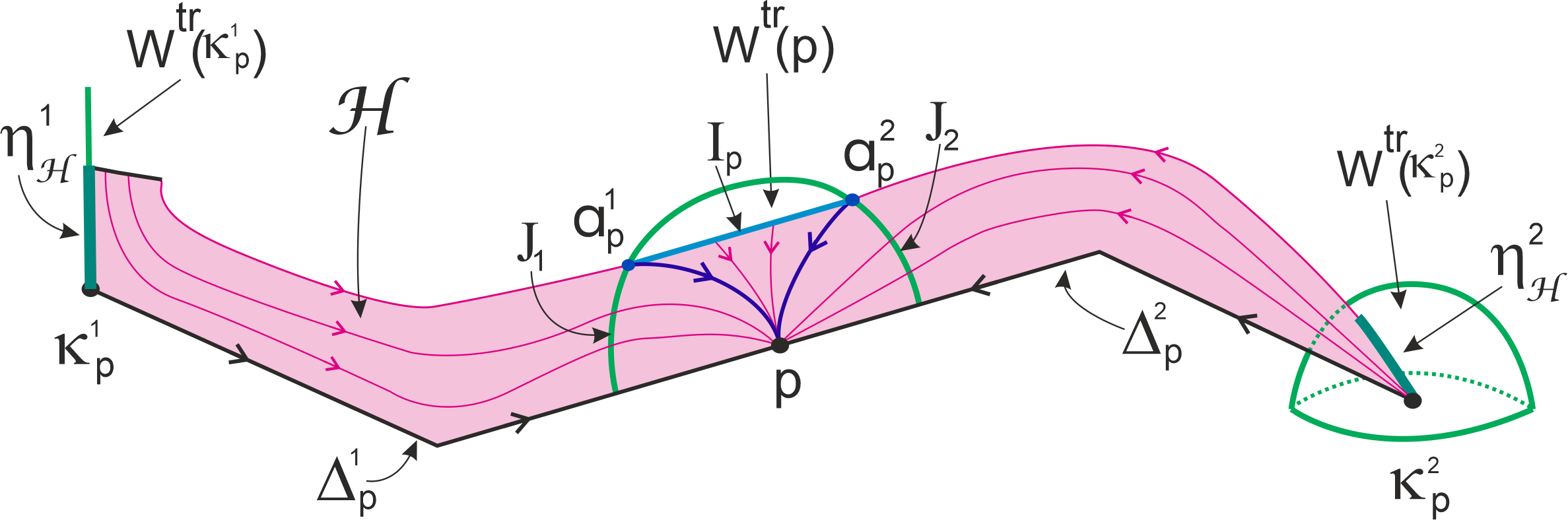}
		\caption{A fixed separating surface and related elements.}
		\label{Fig:SeparantesFijas}
	\end{center}
\end{figure}

Let us consider $\Delta_p^1, \Delta_p^2$ the two paths of edges described in the first part of Proposition~\ref{pro:saturation-marks} that start at points $\kappa^1_p,\kappa^2_p\in N$, respectively, and end at $p$. Denote by $a^1_p,a^2_p$ the extremities of the interval $I_p=T_p\cap\overline{W}$. Taking into account Definition~\ref{def:fitting} (notably the hypothesis that $Fr(F)\cap\overline{W}=T_p\cap\overline{W}$ and the fact that points in $I_p$ are of type t-i relatively to $F$), we have that these extremities $a^i_p$ belong to $\partial W$. Denote by $J_1, J_2 \subset \partial W$ the two open intervals such that $\overline{J_i}\setminus J_i=\{a^i_p, c^i_p\}$ where $c^i_p \in \Delta_p^i\setminus\{p\}$, for $i\in \{1,2\}$ (see Figure \ref{Fig:SeparantesFijas}).

\begin{lemma}\label{lm:FP3}
	 Denote by $\ell^{F,-}_x$ the negative $F$-leaf through a point $x\in F$. The following statements hold for $i\in \{1,2\}$:
	\begin{enumerate}[(a)]
		\item If $\kappa_p^i\in N_3$ (hence $\kappa_p^i\in N_3^u$), then
		$
		\alpha(\ell^{F,-}_x)=\kappa_p^i
		$ for any $x\in J_i$. Moreover, $\alpha(\ell^{F,-}_{a^i_p})=\kappa_p^i$ and $\ell^{F,-}_{a^i_p}\cap Fr(F)$ is a closed interval with extremities $\{b^i_p,a^i_p\}$, where $b_p^i$ is the unique point in $\ell^{F,-}_{a^i_p}\cap\partial T_{\kappa_p^i}$.
		\item If $\kappa_p^i\in N_1$ (hence $\kappa_p^i\in N_1^s$), then $
		\ell^{F,-}_x$ cuts $T_{\kappa_p^i}$ for any $x\in J_i$. Moreover, $\ell^{F,-}_{a^i_p}$ is a closed interval entirely contained in $Fr(F)$ with extremities $\{b^i_p,a^i_p\}$, where $b_p^i$ is the unique point in $\ell^{F,-}_{a^i_p}\cap\partial T_{\kappa_p^i}$. 
	\end{enumerate}	
\end{lemma}
\begin{proof}
	A first statement useful to get both (a) and (b), is the next:
	
\begin{quote}
		{\bf Claim (C).-} If $x\in J_i$ then either $\ell^{F,-}_x$ stays in $U$ and coincides with the whole negative leaf $\ell^-_x$ in $M$, or it escapes from $F$ at some point in $int(T_q)\subset Fr(F)$ for some $q\in N_1^s\cup S_{tr}^s$ . 
\end{quote}

	To prove the claim, assume that $\ell^{F,-}_x$ intersects $Fr(F)$. Let $y$ be the first point (i.e., the closest to $x$) in such intersection. Since $x$ is of type i-i, relatively to $F$, the point $y$ must be of type t-i or e-i. The first case is forbidden by (iv) in Definition~\ref{def:fitting}. We conclude by (iii) in the same definition.
	
	\vspace{.2cm}

	Let us prove the first sentence in (a). Denote $
	B_i=\{x\in J_i\,:\,\alpha(\ell^{F,-}_x)=\kappa_p^i\}$ and let us see that $B_i=J_i$ by showing that $B_i$ is a non-empty open and closed in $J_i$. Since $\kappa_p^i\in|\Delta_p^i|\subset cl_{V}(\Sat_{V}(J_i))$ for some neighborhood $V\subset U$ of $D$, there are points $x\in J_i\cap V$ such that $\kappa_p^i\in cl_V(\ell^{V,-}_x)$. Hence $B_i\ne\emptyset$. Let us show that $B_i$ is open in $J_i$. Take $x\in B_i$. Since $\alpha(\ell^-_x)=\kappa_p^i$ and $\ell_x\cap D=\emptyset$, we have that $\ell^-_x$ cuts the set  $\partial(W^{tr}(\kappa_p^i))\setminus D$, which is an open disc transversal to $\LL$.  By continuity of the flow, it suffices to prove that $\ell^{F,-}_x\subset U$. If it were not the case, by Claim (C) above, there would exist a point $y\in\ell^-_x\in Fr(F)$ which belongs to $int(T_q)$ for some $q\in N^s_1\cup S_{tr}^s$. Hence $\ell^{F,-}_x$ would stop at $y$, contradicting the fact that $x\in B_i$. Finally, let us show that $B_i$ is closed in $J_i$. Suppose that $x\in J_i$ is the limit point of a sequence $\{x_n\}\subset B_i$. Let $y_n$ be the (unique) point where the leaf $\ell^{F,-}_{x_n}$ cuts $\partial (W^{tr}(\kappa_p^i))$. Up to taking a subsequence, we can assume that $\{y_n\}$ converges to some $y\in\partial( W^{tr}(\kappa_i))$. Since $y_n\in\ell^-_{x_n}$ for any $n$, we must have that $y\in\ell^-_{x}$ and hence $\alpha(\ell^-_{x})=\kappa_p^i$. Using Claim (C) and the same reasoning as above, we get  $\ell^-_{x}=\ell^{F,-}_{x}\subset U$ and conclude that $x\in B_i$, as wanted. 
	
	\vspace{.2cm}
	
	To prove the second part of (a), take into account that, since $\alpha(\ell^{-}_{x})=\kappa_p^i$ for any $x\in J_i$, we also have $\alpha(\ell^{-}_{a^i_p})=\kappa_p^i$, by continuity. On the other hand, as $a^i_p$ is a point of type t-i, then $\ell^-_{a^i_p}$ is contained in $Fr(F)$ locally at $a^i_p$ (in fact, in $\overline{Fr(F)^\smile}$). Hence, at some point of $\ell^-_{a^i_p}$ we must switch from type t-t to type e-t or i-t. Let $y$ denote the closest point to $x$ where this change occurs. If $y$ were of type e-t then, according to condition (iii) in Definition~\ref{def:fitting}, we would have $y\in\partial T_q$, with $q\in N^s_1\cup S_{tr}^s$. In both cases, the negative $F$-leaf $\ell^{F,-}_{a^i_p}$ would stop at $y$ so that it cannot reach $\kappa_p^i$. This would also happen for $\ell^{F,-}_x$, with $x\in J_i$ sufficiently near to $a^i_p$. Contradiction. Thus, necessarily $y$ is of type i-t. By conditions (iii) and (iv) of Definition~\ref{def:fitting}, we must have $y\in\partial T_q$ for some $q\in N_3^u$. We conclude that  $q=\kappa_p^i$ and we are done (with the extremity $b^i_p$ in statement (a) equal to the point $y$). 
	
	\vspace{.2cm}
	
	Let us see the first sentence in (b). Denote by $
	B_i=\{x\in J_i\,:\,\ell^{F,-}_x\mbox{ cuts }T_{\kappa_p^i}\}
	$ and let us show again that $B_i$ is non-empty open and closed in $J_i$. By means of Hartman-Grobman's Theorem, if $V_i$ is a sufficiently small neighborhood of $\kappa_p^i$ and $z\in V_i\setminus D$ (i.e., not in the unstable manifold at $\kappa_p^i$ ), then $\ell_z$ cuts $T_{\kappa_p^i}$ (which is a transversal disc attached at some point in the stable manifold $W^s_{\kappa_p^i}$). Since $\kappa_ p^i\in|\Delta_p^i|\subset cl_V(\Sat_V(J_i))$, if $V$ is a small neighborhood of $D$ in $M$, we conclude that $B_i\ne\emptyset$. To show that $B_i$ is closed in $J_i$, let $x\in J_i$ be the limit point of a sequence $\{x_n\}\subset B_i$ so that $\ell^{F,-}_{x_n}$ cuts $T_{\kappa_p^i}$ (necessarily in a unique point, since the points of $T_{\kappa_p^i}$ are of type e-t or e-i relatively to $F$). By continuity of the flow, the whole negative leaf  $\ell^-_x$ must cut also $T_{\kappa_{p}^i}$. We need to show that the point in $\ell^-_x\cap T_{\kappa_p^i}$ closest to $x$ belongs to the negative $F$-leaf $\ell^{F,-}_x$. If not, using Claim (C), $\ell^{F,-}_x\cap Fr(F)$ would be a point in $int(T_q)$ with some $q\ne\kappa_p^i$. By continuity, the negative leaf $\ell^{-}_{x_n}$ also cuts $int(T_q)$ if $n$ is big enough, so that $\ell^{F,-}_{x_n}$ stops at $int(T_q)$ before reaching $T_{\kappa_p^i}$. Contradiction. Finally,  if $x\in B_i$ and $y\in\ell^{F,-}_x\cap T_{\kappa_p^i}$ then, using again Claim (C) and similar arguments, we must have that $y\in int(T_{\kappa_p^i})$. This shows that $B_i$ is open in $J_i$ and we are done.
	
	\vspace{.2cm}
	
	Finally, we prove the last part in (b). By continuity, the leaf $\ell^{-}_{a^i_p}$ cuts $T_{\kappa_ p^i}$. Let $y$ be the closest point to $x$ with this property. As done in the proof of the second sentence in (a), we see that the type, relatively to $F$, remains constantly equal to t-t along all the points in the open interval in $\ell^-_{a^{i}_p}$ between $y$ and $x$. On the other hand, $\ell^-_{a^i_p}$ escapes from $F$ at the point $y$, since $y$ is of type e-t. Consequently, $\ell^{F,-}_{a^i_p}$ is equal to the closure of the interval in the leaf $\ell^-_{a^i_p}$ between $b^i_p:=y$ and $a^i_p$. This ends the proof of (b). 
\end{proof}

In the next result we describe the main properties of  fixed separating surfaces.
\begin{proposition}\label{pro:fixed-separating}
Assume that $W=W^{tr}(p)=W^s(p)$ is the stable manifold at $p$. Then the fixed separating surface $\HH$ satisfies the following properties:
\begin{enumerate}[(FS1)]
	\item It is a smooth analytic two-dimensional connected submanifold of $U\setminus D$.
	\item For any $x\in\HH$, the positive leaf $\ell^+_x$ of $\LL$ in $M$ through $x$ is contained in $\HH$ and satisfies $\omega(\ell^+_x)=p$.
	\item  The closure $cl_U(\HH)=cl(\HH)\cap U$ of $\HH$ in $U$ is a topological surface with connected boundary. More precisely, there are $U$-leaves $\eta^1_{\HH},\eta^2_\HH$ in $U\setminus D$ with $cl_U(\eta^i_{\HH})\setminus\eta^i_{\HH}=
	\{\kappa^i_{p}\}$ for $i=1,2$, such that the boundary $\partial\HH$ of $cl_U(\HH)$ is given by
	\begin{equation}\label{eq:cl-fixed-surface}
		\partial\HH=cl_U(\HH)\setminus\HH=
		|\Delta_p^1|\cup|\Delta_p^2|
		\cup\eta^1_{\HH}\cup\eta^2_\HH.
	\end{equation} 
	Moreover, in case $\kappa_p^i\in N_1$, then  $\eta^i_{\HH}=W^{tr}(\kappa^i_p)\cap U$ (see Figure \ref{Fig:SeparantesFijas}).

	\item For each $i\in\{1,2\}$, the intersection $\ell^U_{a^i_p}:=\ell_{a^i_p}\cap U$ is an $U$-leaf  contained in $\HH$. Moreover, the set $\HH\setminus(\ell^U_{a^1_p}\cup\ell^U_{a^2_p})$ has three connected components, all of them saturated by $\LL$ in $U$, which can be indexed as $\HH^0,\HH^1,\HH^2$ in such a way that they satisfy: 
	\begin{enumerate}[(a)]
		\item If $x\in\HH^0$, then the negative leaf $\ell^-_x$ scapes from $U$ at some point of $I_p$.
		\item For $i\in\{1,2\}$, if $\kappa^i_p\in N^u_3$ then, for any $x\in\HH^i$, the negative leaf $\ell^-_x$ is entirely contained in $\HH^i$ and satisfies $\alpha(\ell^-_x)=\kappa_p^i$.
		\item For $i\in\{1,2\}$, if $\kappa^i_p\in N^u_1$ then, for any $x\in\HH^i$, the negative leaf $\ell^-_x$ scapes from $U$ at some point of the transversal disc $T_{\kappa_p^i}$. 
	\end{enumerate} 
\end{enumerate}
\end{proposition}
\begin{proof}

We have $\HH=\Sat_U(W\cap U\setminus D)$, where we have denoted $W:=W^{tr}(p)$. Recall that $\overline{W}$ is a closed disc and that $\partial\overline{W}\cap D=W\cap D$ is a non-trivial interval. 

\vspace{.2cm}

\paragraph{\textit{Proof of (FS1)}}

Since $W\setminus D$ is an analytic smooth surface of $int(M)$, the set $(W\cap U)\setminus D$ is an analytic smooth surface of $U\setminus D$. By the analyticity of the flow, the same is true for $\HH$.
To show that $\HH$ is connected, it suffices to prove that $W\cap U$ is. Given that $\LL$ is transversal to the boundary of $\overline{W}$ along $\partial \overline{W}\setminus D$ and that $W$ is the stable manifold at $p$,  the following hold for any $a\in\overline{W}$:
\begin{itemize}
	\item The $\overline{W}$-leaf through $a$ is just $\ell_a\cap\overline{W}$.
	\item The positive leaf $\ell^+_a$ is contained in $\overline{W}$ and satisfies $\omega(\ell^+_a)=\{p\}$. 
	\item When $a\not\in D$, the negative $\overline{W}$-leaf $\ell^-_a\cap\overline{W}$ cuts just once the boundary $\partial\overline{W}$.
\end{itemize}

Notice that the interval $I_p$ is contained in $\overline{W}\setminus D$ and its extremities $a^1_p,a^2_p$ belong to $\partial\overline{W}$ (cf. Definition~\ref{def:fitting}). The map $I_p\to\partial W$ defined by $a\mapsto\ell^-_a\cap\partial W$ is injective and continuous, so that its image is a closed interval $J_p$ in $\partial\overline{W}$. Let
$
K:=\bigcup_{a\in I_p}\ell^-_a\cap\overline{W},
$
a closed semi-analytic set, negatively saturated inside $\overline{W}$. 
\begin{figure}[h]
	\begin{center}
		\includegraphics[scale=0.6]{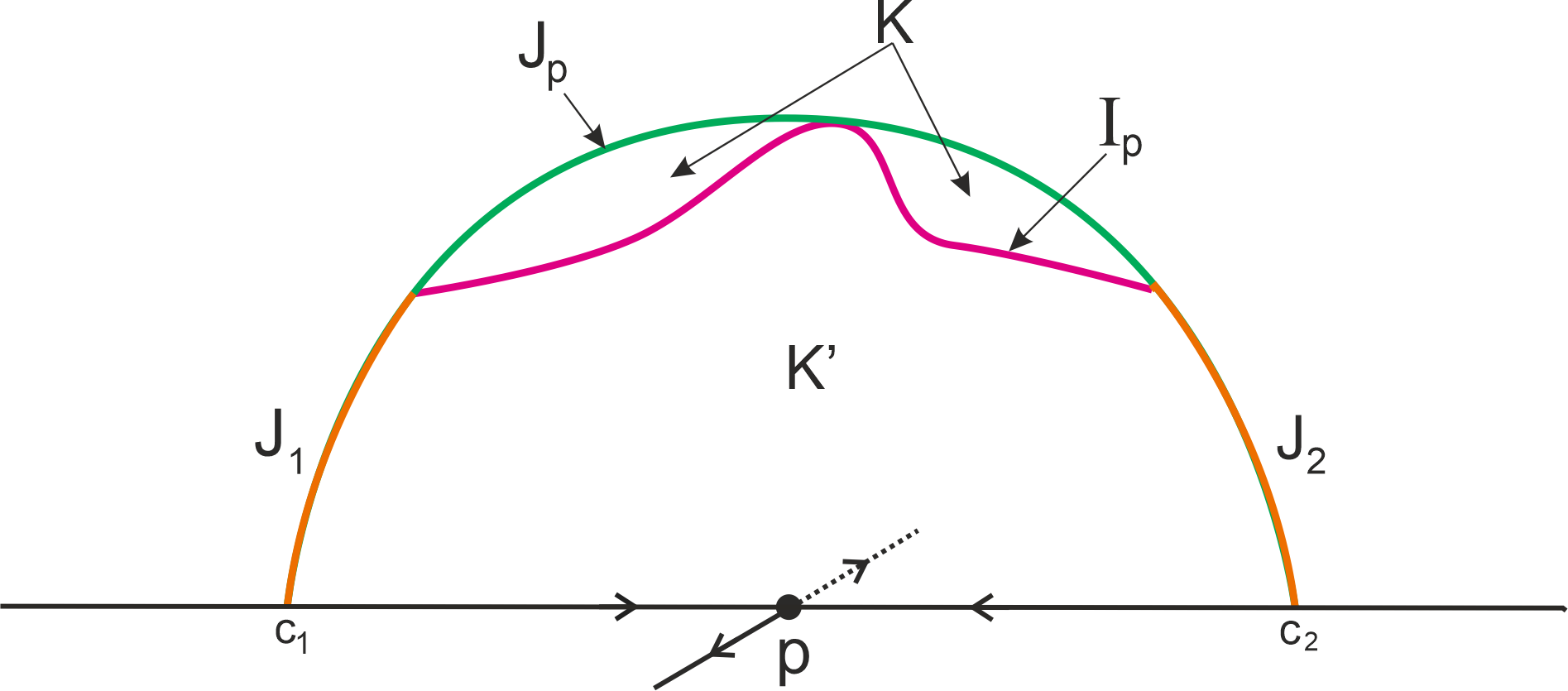}
		\caption{Distinc elements in $W$.}
		\label{Fig:PlainDisc}
	\end{center}
\end{figure}

 Using the above properties, it can be shown without difficulty that $K':=\overline{W\setminus K}$ is a closed disc and that its boundary is equal to (see Figure \ref{Fig:PlainDisc}).

\begin{equation}\label{eq:frontierK-prima}
	\partial K'=I_p\cup(\overline{W}\cap D)\cup J_1\cup J_2,
\end{equation}
where the intervals $J_1,J_2$ are those considered in Lemma~\ref{lm:FP3}. 
Let us show that $int(K')=W\cap U\setminus D$, which finishes the proof that $W\cap U\setminus D$ is connected.
Let $x\in int(K')$. Hence, clearly $x\not\in D$. Also, as $x\in\overline{W}$, if $x$ were in $\partial\overline{W}$, using (\ref{eq:frontierK-prima}), we would have that $x\in J_p\subset K$, which is impossible. Thus $x\in W$. On the other hand, since $\omega(\ell^+_x)=\{p\}$, if $x\not\in U$ then $\ell^+_x$ cuts the frontier of $U$. By the definition of fitting domain, a point in $\overline{W}\cap Fr(U)$ is either of type t-t, e-i or t-i, relatively to $F$. In the last two cases, such a point belongs to $I_p$, but this is impossible since $\ell^+_x\cap I_p=\emptyset$ (otherwise $x\in K$). Thus, we must have  that all points of $\ell^+_x$ are of type t-t, and hence $\ell^+_x\subset Fr(U)$, which is again impossible because the accumulation point $p$ does not belong to $Fr(U)$. We deduce the inclusion $int(K')\subset W\cap U\setminus D$. Take now a point $x\in W\cap U\setminus D$. Then $\ell^+_x$ cannot cut $I_p$ (since $I_p$ is composed of points of type e-i or t-i), so that $x\not\in K$ and hence $x\in K'$. Also, using (\ref{eq:frontierK-prima}), and taking into account that $\partial K'\cap(W\setminus D)=I_p\cap W$ and $I_p\subset K$, we conclude that $x\not\in\partial K'$ and hence $x\in int(K')$.

\vspace{.2cm}

\paragraph{\textit{Proof of (FS2)}}
This is quite obvious: such a property holds for any $x\in W\cap U$ and, by definition, $\HH=\Sat_{U}(W\cap U\setminus D)$.

\vspace{.2cm}

\paragraph{\textit{Proof of (FS3)}}

By virtue of Proposition~\ref{pro:saturation-marks}, there is a neighborhood $V\subset U$ of $D$ satisfying $cl_{V}(\HH\cap V)\cap D=|\Delta_p^1|\cup|\Delta_p^2|$. Hence $|\Delta_p^1|\cup|\Delta_p^2|\subset cl_U(\HH)\setminus\HH$. This is a first step towards (\ref{eq:cl-fixed-surface}).

We define the leaves $\eta^i_\HH$ using the notations and cases appearing in the statement of Lemma~\ref{lm:FP3} (see again Figure \ref{Fig:SeparantesFijas}). 
\begin{itemize}
	\item If $\kappa_p^i\in N_3^u$ (case (a)), we put 
	$$
	\eta^i_\HH:=\ell^{F,-}_{b^i_p}\setminus\{b^i_p\}.
	$$
	\item If $\kappa_p^i\in N_1^s$ (case (b)), we put 
	$$
	\eta^i_\HH:=W^{tr}(\kappa_p^i)\cap U.
	$$
	
\end{itemize}

\noindent Notice that, in both cases, $\eta^i_\HH$ is an $U$-leaf. Following the proof of Lemma~\ref{lm:FP3},  in case (a) we can assert that $\eta^i_{\HH}=\ell^-_{a^i_p}\cap U$. As a consequence, since $\ell^-_x\subset U$ for any $x\in J_i$, we have that $\eta^i_\HH$ is contained in $cl_U(\HH)$. If we are in case (b), we have also the inclusion  $\eta^i_\HH\subset cl_U(\HH)$ as a consequence of Hartman-Grobman's Theorem, since $\kappa^i_p\in cl_U(\HH)$. 
To end the proof of (\ref{eq:cl-fixed-surface}) it remains to show the inclusion
$$
cl_U(\HH)\setminus \HH\subset \eta_\HH^1\cup\eta_\HH^2\cup|\Delta_p^1|\cup|\Delta_p^2|.
$$
Let $z\in cl_U(\HH)\setminus\HH$. If $z\in D$ then, by Proposition~\ref{pro:saturation-marks}, we have that $z\in|\Delta_p^1|\cup|\Delta_p^2|$ and we are done. Otherwise, supose that $z\not\in D$ and consider a sequence $\{z_n\}\subset \HH$ converging to $z$. By definition of $\HH$, for each $n$, there is a point $x_n\in W\cap U\setminus D$ such that $z_n\in\ell_{x_n}$. Let $y_n$ be the unique point where the leaf $\ell_{x_n}(=\ell_{z_n})$ cuts $\partial W$ (notice that the boundary of $W$ is transversal to $\LL$ along $\partial W\setminus D$ and that $x_n\not\in D$). Let $y\in\partial W$ be an accumulation point of the sequence $\{y_n\}$. We have two cases:
\begin{enumerate}
	\item If $y\not\in D$, one can see by continuity of the flow, that $z\in\ell_y$. Since $z\in U\setminus \HH$, we must have $y\not\in I_p$ (otherwise $\ell_y\cap U=\ell^{U,+}_y\setminus\{y\}\subset int(W)$ as we have seen in the proof of (FS1). Analogously, by Lemma~\ref{lm:FP3}, we have $y\not\in J_1\cup J_2$. Thus $y=a^i_p$ for some $i\in\{1,2\}$. We deduce, using again Lemma~\ref{lm:FP3}, that $\kappa_p^i\in N_3^u$ (case (a)) and that $z\in\eta_\HH^{i}$ (since $\eta_\HH^i=\ell^-_{a_p^i}\cap U$ in this case).
	\item If $y\in D$, since $cl_V(\HH\cap V))\cap D=|\Delta_p^1|\cup|\Delta_p^2|$ for a  small enough neighborhood $V$ of $D$, there is some $i\in\{1,2\}$ and an open neighborhood $V_i$ of $\kappa_p^i$ such that, for any $n$ sufficiently large, the leaf $\ell_{y_n}=\ell_{x_n}$ intersects $V_i$ and  the family of local leaves $\{\ell_{y_n}\cap V_i\}_n$ accumulates along the set $|\Delta_p^i|\cap V_i$. We must have then that $\kappa_p^i\in N_1^s$ (case (b)): in case $\kappa_p^i$ were in $N_3^u$ we would have that $\alpha(\ell_{y_n})=\kappa_p^i$, for any $n$ large enough, so that the entire leaf $\ell_{y_n}=\ell_{x_n}$ is included in $U$ (by Lemma~\ref{lm:FP3}, (a)) and the family $\{\ell_{y_n}\}_n$ accumulates only along $D$, contradicting the fact that $z=\lim z_n$ with $z_n\in\ell_{y_n}$ and $z\not\in D$. 
	Also, by Lemma~\ref{lm:FP3}, (b), we may impose that the sequence $\{z_n\}$ is contained in $V_i$, if we assume that $V_i$ contains $T_{\kappa_p^i}$. Otherwise, $\ell_{y_n}\cap (U\setminus V_i)$ would contain an $(U\setminus V_i)$-leaf that accumulates along $D$ when $n\to\infty$, and we would have again the contradiction $z=\lim z_n\in D$. Given that the local leaves $\ell_{y_n}\cap V_i$ accumulate to
	the subset $|\Delta_p^i|\cap V_i$, contained in the local unstable manifold $W^u_{\kappa_p^i}=D\cap V_i$ at $\kappa_p^i$ we must have, again by the Hartman-Grobman's Theorem, that they accumulate also to the stable manifold $W^s_{\kappa_p^i}=\eta_\HH^i\cup\{\kappa_p^i\}$ (recall that we are in case (b)). We conclude that $z\in\eta_\HH^i$, as wanted.
\end{enumerate}

\vspace{.2cm}

\paragraph{\textit{Proof of (FS4)}}

The properties stated in this item are essentially already proved in (FS1)-(FS3). We need to observe that the leaf $\ell^U_{a^i_p}$ has extremities $\{a^i_p,p\}$ and that the union of $\ell^U_{a^1_p}$ and $\ell^U_{a^2_p}$ separates $int(K')=W\cap U\setminus D$ (cf. proof of (FS1)) into three connected components $K'_0,K'_1,K'_2$ whose closures cut the transversal boundary $\partial\overline{W}\setminus D$ along $I_p,\overline{J_1},\overline{J_2}$, respectively. Since any $U$-leaf inside $\HH$ cuts $int(K')$, we get the decomposition $\HH\setminus(\ell^U_{a^1_p}\cup\ell^U_{a^2_p})=\HH^0\cup\HH^1\cup\HH^2$, where $\HH^j=\Sat_{U}(K'_j)$, $j=0,1,2$. On the other hand, since the leaf through any point in $\HH$ cuts $\partial K'\setminus D$, properties (a), (b), (c) of (FS4) are consequences of (FS1) along with items (a) and (b) of Lemma~\ref{lm:FP3}.

This concludes the proof of Proposition~\ref{pro:fixed-separating}.
\end{proof}

\subsection{Mobile separating surfaces properties}\label{subsec:prop-mobile-separating}
Let $\GG:=\GG^{i,\epsilon}_p$ be a mobile separating surface generated at some $p\in S_{tr}$. We gather the main properties of $\GG$ in the main result of this subsection. As before, we assume  that $W:=W^{tr}(p)$ is the stable manifold at $p$. The corresponding properties in case that it is the unstable manifold hold analogously.

\begin{figure}[h]
	\begin{center}
		\includegraphics[scale=0.6]{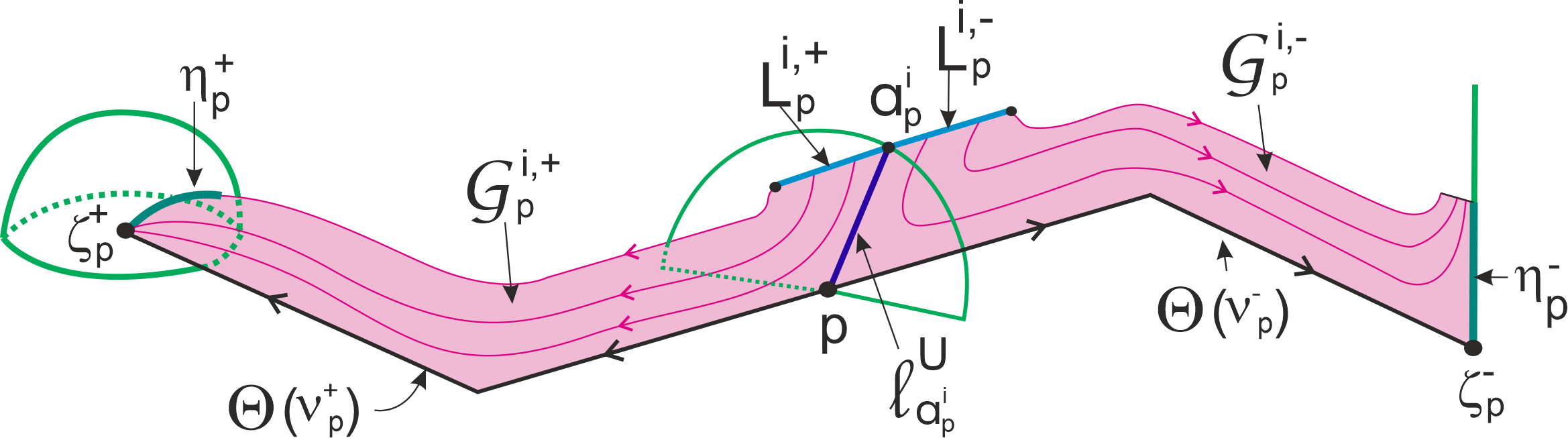}
		\caption{A couple of separating mobile surfaces and related elements.}
		\label{Fig:SeparantesMoviles}
	\end{center}
\end{figure}

With the notations in  Definition~\ref{def:fitting}, recall that $\GG=\Sat_F(L^{i,\epsilon}_p)\cap U$, where $L^{i,\epsilon}_p$ is one of the connected components of $L^i_p\setminus\{a^i_p\}$ and $L^i_p$ is the interval in $\partial T_p$ of points of type t-i. 
Denote also by $\nu^\epsilon_p$ the local side of $W$ at $p$ satisfying $L^{i,\epsilon}_p\subset\nu^\epsilon_p$ and let $\Theta(\nu^\epsilon_p)$ be the path of edges defined in Proposition~\ref{pro:saturation-marks}. Notice that $\Theta(\nu^\epsilon_p)$ goes from $p$ to some vertex in $N^s$,  in what follows denoted by $\zeta^\epsilon_p$ (see Figure \ref{Fig:SeparantesMoviles}). Recall that it does not depend on $i\in\{1,2\}$, just  on $p$ and $\epsilon$.

\begin{proposition}\label{pro:mobile-separating}
Assume that $W^{tr}(p)=W^s(p)$. Then the mobile separating surface $\GG$ satisfies the following properties:
\begin{enumerate}[(MS1)]
 \item It is a two-dimensional smooth analytic connected submanifold of $U\setminus D$.
 \item For any $x\in\GG$, the negative $U$-leaf through $x$ scapes from $U$ at a point in $L^{i,\epsilon}_p$.
 \item The closure $cl_U(\GG)=cl(\GG)\cap U$ is a topological surface with boundary. More precisely, there exists a leaf $\eta_{\GG}$  in $U\setminus D$ satisfying $cl_U(\eta_{\GG})\setminus\eta_{\GG}
 =\{\zeta^\epsilon_{p}\}$ such that the boundary $\partial\GG$ of $cl_U(\GG)$ is given by
 \begin{equation}\label{eq:cl-mobile-surface}
 	\partial\GG=cl_U(\GG)\setminus\GG=
 	|\Theta(\nu^\epsilon_p)|\cup
 	\ell^U_{a^i_p}\cup\eta_\GG.
 \end{equation}
Moreover, if $\zeta^\epsilon_p\in N_1$, then we have $\eta_\GG=W^1_{\zeta^\epsilon_p}\cap U$. 
\item The asymptotic behavior of the positive leaves in $\GG$ is described as follows:
\begin{enumerate}[(i)]
  \item If $\zeta^\epsilon_p\in N_3$ (necessarily in $N^s_3$) then, for any $x\in\GG$, the positive leaf $\ell^+_x$ is entirely contained in $\GG$ and satisfies   $\omega(\ell^+_x)=\zeta^\epsilon_p$. 
  \item If $\zeta^\epsilon_p\in N_1$ (necessarily in $N^u_1$) then, for any $x\in\GG$, the positive leaf $\ell^+_x$  scapes from $U$  at some point of the transversal disc $T_{\zeta^\epsilon_p}$. 
\end{enumerate} 
\end{enumerate}
\end{proposition}
\begin{proof}
We proceed as in the proof of  Proposition~\ref{pro:fixed-separating}. Let us just indicate the main steps. Consider first $\{a^i_p,d=d^{i,\epsilon}_p\}$, the set of extremities of the interval $L^{i,\epsilon}_p$. Recall that $d$ is a point of type t-t relatively to $F$, whereas any other point of $L^{i,\epsilon}_p$ is of type t-i (cf. Definition~\ref{def:fitting}). Analogously to the proof in Lemma~\ref{lm:FP3}, we prove that, for any $x\in L^{i,\epsilon}_p\cup\{d\}$, the positive $F$-leaf $\ell_{x}^{F,+}$ satisfies:
\begin{enumerate}[(a)]
	\item If $\zeta^\epsilon_p\in N^s_3$, then $\ell_{x}^{F,+}=\ell^+_x$ and $\omega(\ell^+_x)=\zeta^\epsilon_p$. Moreover, $\ell^+_x\cap U=\ell^+_x\setminus\{x\}$ if $x\ne d$. This last property shows (i) of item (MS4).
	
	\item  If $\zeta^\epsilon_p\in N^s_1$, then $\ell_{x}^{F,+}=\ell^+_x\cap F$ and $\ell_{x}^{F,+}$ is a closed interval with an extremity in $x$ and the other one in $T_{\zeta^\epsilon_p}$. Moreover, if $x=d$, then $\ell_{d}^{F,+}$ is contained in $Fr(F)$; otherwise, $\ell_{x}^{F,+}\cap U=int(\ell_{x}^{F,+})$. This proves (ii) of item (MS4).
\end{enumerate}

It is clear that, in the previous situations (a) and (b), the mobile separating surface $\GG$ may be defined also as $\GG=\Sat_U(L'^{i,\epsilon}_p)$, where $L'^{i,\epsilon}_p$ is the image of $L^{i,\epsilon}_p$ by the flow of (a generating vector field of) $\LL$ at a small time $\delta>0$. Since $L'^{i,\epsilon}_p$ is also an open interval, we obtain that $\GG$ is connected. The rest of properties stated in (MS1) are straightforward. Item (MS2) is also consequence of (a) and (b). Finally, for item (MS3), we define the $U$-leaf $\eta_\GG$ as follows (see Figure \ref{Fig:SeparantesMoviles}):
\begin{itemize}
	\item If $\zeta^\epsilon_p\in N_3$ (case (a)), we put $\eta_\GG:=\ell^+_{c}\cap U$ where $c=c^{i,\epsilon}_p$ is the unique point where $\ell^{+}_{d}$ intersects $W^{tr}(\zeta^\epsilon_p)$ (notice that $c$ belongs to $\partial T_{\zeta^\epsilon_p}$).
	\item If $\zeta^\epsilon_p\in N_1$ (case (b)), we put $\eta_{\GG}:=W^{tr}(\zeta^\epsilon_p)\cap U$.
\end{itemize}

With this definition, we prove equation (\ref{eq:cl-mobile-surface}): the inclusion $\Theta(\nu^\epsilon)\cup\eta_\GG\cup\ell^U_{a^i_p}\subset cl_U(\GG)$ is consequence of Proposition~\ref{pro:saturation-marks} and the construction of $\eta_\GG$ and $\ell^U_{a^i_p}$; the other inclusion is shown using similar arguments as those used to get equation (\ref{eq:cl-fixed-surface}) at the end of the proof of (FS3).
\end{proof}

\subsection{The strata. End of the proof of Theorem~\ref{th:afterblowup-intro}}
Recall that we had considered a neighborhood basis $\{F_n\}_n$ of $D$, where each $F_n$ is a fitting domain. To finish the proof of  Theorem~\ref{th:afterblowup-intro}, we define an analytic stratification $\mathcal{S}_n$ of $U_n:=int(F_n)$ into dynamically elementary pieces satisfying all the required properties. To this end, we use the notations in Proposition~\ref{pro:fixed-separating} and Proposition~\ref{pro:mobile-separating} concerning the fixed and mobile separating surfaces of each element $F_n$ in that basis. To lighten the notation, we remove again the index $n$ from $F_n$, $U_n$ and the stratification $\SC_n$, as well as for the elements that are associated to the fitting domain $F_n$, for instance the transversal disc $T_q$, etc. Nevertheless, we must keep in mind that we are, in fact, defining a sequence of stratifications in the different open sets $U_n$, for $n\in\N$, so that expressions such as ``... does not depend on $n$'' makes sense (cf. items (i) and (iii) in Theorem~\ref{th:afterblowup-intro}).

\begin{figure}[h]
	\begin{center}
		\includegraphics[scale=0.6]{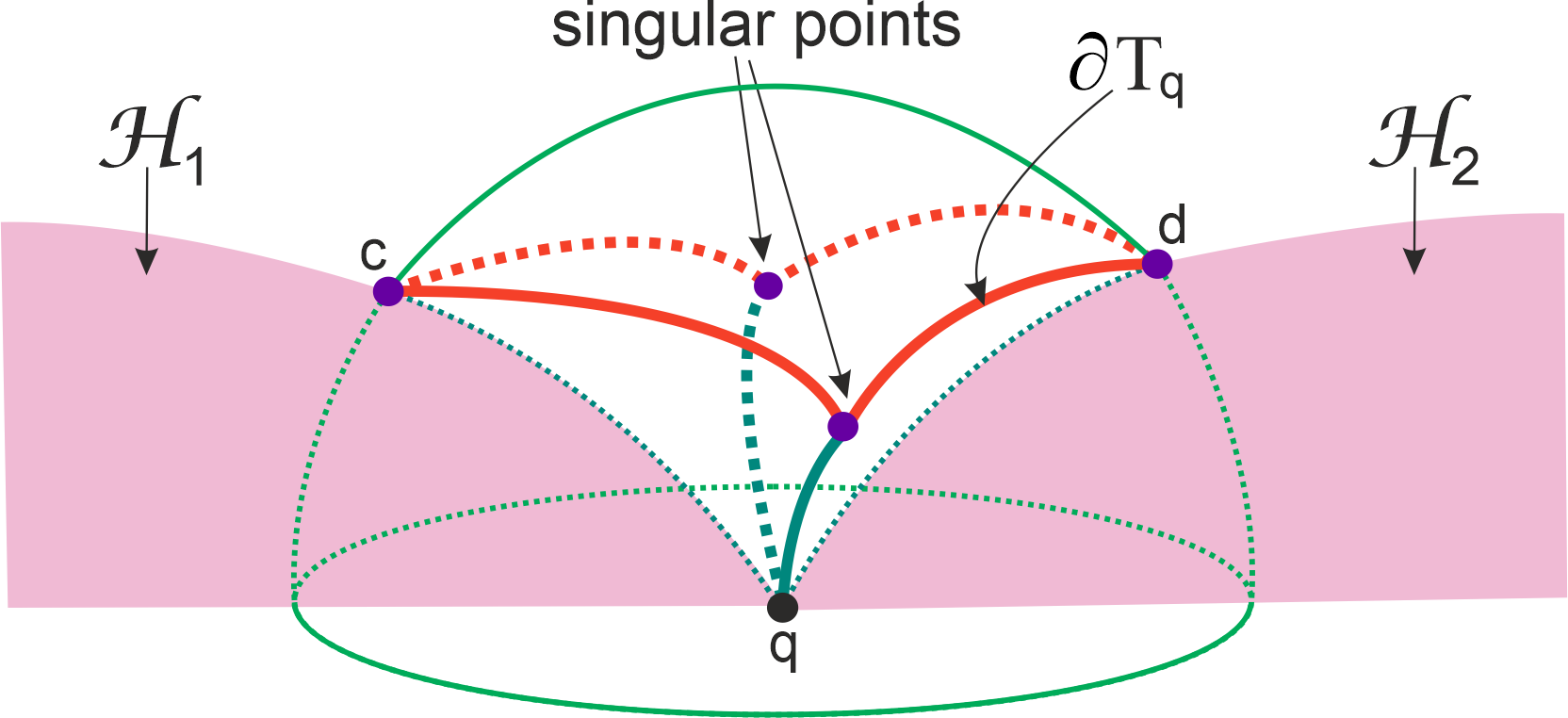}
		\caption{Local parabolic strata at a point $q\in N_3$. The separa\-ting surfaces $\mathcal{H}_1, \mathcal{H}_2$ cut $\partial T_q$ at points $c,d\in\Upsilon_q$.}
		\label{Fig:RefinementStrata}
	\end{center}
\end{figure}

\strut
With these conventions, the elements of the stratification $\SC$ of $U=int(F)$ are the following ones:

\begin{enumerate}
	\item[(d0)] {\em Zero-dimensional} strata. They are exactly the singular points of $\LL$ or, that is, the vertices of the associated graph $\Omega$. 
	\item[(d1)] {\em One-dimensional} strata. They are the sets belonging to one of the following families:
	\begin{enumerate}
		\item[(1.a)] Any edge of the graph $\Omega$.
		\item[(1.b)] Any of the $U$-leaves
		$\ell^U_{a^1_p}$ or $\ell^U_{a^2_p}$, where $p \in S_{tr}$.
		\item[(1.c)] Any of the $U$-leaves $\eta^1_\HH$ or $\eta^2_\HH$, where $\HH$ runs in the set of fixed separating surfaces of $F$ (cf. Proposition~\ref{pro:fixed-separating}, (FS3)).
		\item[(1.d)] Any $U$-leaf $\eta_\GG$, where $\GG$ runs in the set of mobile separating surfaces of $F$ (cf. Proposition~\ref{pro:mobile-separating}, (MS3)).
		\item[(1.e)] Any $U$-leaf of the form $U\cap W^1_q$, where $q\in N_1$, not belonging to the previous families (1.c) or (1.d).
		\item[(1.f)] Any $U$-leaf of the form $\ell_x\cap U$, where $x\in \Sing(\partial T_q)$ and $q\in N_3$ (and which is not already in the families (1.c) or (1.d)). Here, $\Sing(\partial T_q)$ denotes the set of points in $\partial T^n_q$ where this curve is not smooth analytic. (Notice that, since $\partial T_q$ is compact and semi-analytic, the set $\Sing(\partial T_q)$ is finite). 

	\end{enumerate}
	\item[(d2)] {\em Two-dimensional} strata.  They are the sets in one of the following families:
	\begin{enumerate}
		\item[(2.a)] Any face of the graph $\Omega$.
		\item[(2.b)] Any of the components $\HH^0,\HH^1,\HH^2$ described in Proposition~\ref{pro:fixed-separating}, (FS4), where $\HH$ runs in the set of fixed separating surfaces.
		\item[(2.c)] Any mobile separating surface 
		$\GG=\GG^{i,\epsilon}_p$ for $i\in\{1,2\}$, $\epsilon\in\{+,-\}$, where $p \in S_{tr}$.
		\item[(2.d)] 
		Given $q\in N_3$, denote by $\Upsilon_q\subset\partial T_{q}$ the set of points where $\partial T_q$ cuts the closure of some separating surface of $F$ (equal to the set of points in $\partial T_q$ that are extremities of strata in the families (1.c) or (1.d)). Then, any set of the form $\Sat_F(B)\cap U$, where $B$ is a connected component of $\partial T_q\setminus(\Upsilon_q\cup\Sing(\partial T_q))$ and $q \in N_3$, is a two-dimensional stratum of $\SC$ (see Figure~\ref{Fig:RefinementStrata}).
	\end{enumerate}
	\item[(d3)] {\em Three-dimensional} strata. 
	Any connected component of $U\setminus\mathcal{T}^2$, where 
	 $\mathcal{T}^2$ is the union of the subsets of $U$ defined in items (d0), (d1) and (d2). 
\end{enumerate}

From their definition, the family $\SC=\{A_j\}_{j\in J}$ of all these sets forms a partition of $U$. With the help of Proposition~\ref{pro:fixed-separating} and Proposition~\ref{pro:mobile-separating}, we can see that $\SC$ is in fact an analytic $\LL$-stratification, where the strata of dimension $k\in\{0,1,2,3\}$ are those considered in step (dk). More in precise, let $A=A_j$ be one element of this partition:

\begin{itemize}
	\item If $\dim(A)=0$, then $A$ is  clearly saturated by $\LL$ and an elementary dynamical piece.
	\item Suppose that $\dim A=1$. If $A$ is an edge of $\Omega$ (type (1.a)), then $A\subset D$ is a leaf (in $U$, but also in $M$) and $\overline{A}\setminus A$ is composed of the two different points $\alpha(A),\omega(A)$ (thus $A$ is an elliptic elementary dynamical piece). In the other cases (1.b), (1.c), (1.d) and (1.e), $A$ is an $U$-leaf contained in $A\subset U\setminus D$ and $cl_U(A)\setminus A$ is a single point, either equal to $\omega(A)$ or $\alpha(A)$, by  construction. Thus $A$ is a parabolic elementary dynamical piece. In all the cases, by analyticity of $\LL$, $A$ is an analytic submanifold of $U$.
	\item If $A\subset\Omega$ is a face of $\Omega$ (type (2.a)), then $A$ is $\LL$-saturated and, being an open set of $D$, it is an analytic submanifold of $U$. On the other hand, thanks to the acyclicity of the graph (see \cite[Lemma 4]{Alo-San1}), we have that $cl(A)\setminus A$ is equal to the union of supports of two paths of edges of $\Omega$ going both from $q_1$ to $q_2$, where $q_1,q_2$ are two $D$-node points; thus the frontier of $A$ is a union of strata of $\SC$ of smaller dimension. Moreover, any leaf $\ell$ contained in $A$ satisfies $\alpha(\ell)=q_1$ and $\omega(\ell)=q_2$. Thus, $A$ is an elliptic elementary dynamical piece.
	\item If $A=\HH^0$ for some fixed separating mark $\HH$ generated at $p\in S_{tr}$ (type (2.b)) then $A$ is saturated in $U$, since it is a connected component of a difference of two saturated sets in $U$. Being an open subset of $\HH$, it is an analytic submanifold of $U$ (by Proposition~\ref{pro:fixed-separating}, (FS1)). On the other hand, we have $cl_U(A)\setminus A=\ell^U_{a^1_p}\cup\ell^U_{a^2_p}\cup\{p\}$, again a union of strata of smaller dimension. Also, by Proposition~\ref{pro:fixed-separating} ((FS2) and (FS4)), $A$ is a parabolic elementary dynamical piece with $\{\alpha(A),\omega(A)\}=\{\emptyset,p\}$. 
	\item If $A=\HH^i$ with $i\in\{1,2\}$, where $\HH=\HH_p$ is the fixed separating surface generated at $p$ (type (2.b)), then, as above, $A$ is saturated in $U$ and an analytic submanifold of $U$. Also, by construction and (\ref{eq:cl-fixed-surface}), we have that $cl_U(A)\setminus A=\ell^U_{a^i_p}\cup|\Delta_p^i|\cup\eta_\HH^i$, a union of strata of smaller dimension. In light of Proposition~\ref{pro:fixed-separating} ((FS2) and (FS4)), we deduce that $A$ is an elementary dynamical piece. Moreover, $A$ is elliptic when $\kappa^i_p\in N_3$ (with $\{\alpha(A),\omega(A)\}=\{p,\kappa_p^i\}$), or parabolic when $\kappa^i_p\in N_1$ (with $\{\alpha(A),\omega(A)\}=\{\emptyset,p\}$).
	\item  If $A=\GG$, where $\GG=\GG^{i,\epsilon}_p$ is a mobile separating surface generated at $p$ (type (2.c)), then $A$ is saturated in $U$ and an analytic submanifold of $U$, by Proposition~\ref{pro:mobile-separating}.  Also, using equation (\ref{eq:cl-mobile-surface}), we have $cl_U(A)\setminus A=\ell^{U}_{a^i_p}\cup|\Theta(\nu^\epsilon_p)|\cup\eta_\GG$, a union of strata of smaller dimension. By Proposition~\ref{pro:mobile-separating} ((MS2) and (MS4)), we have that $A$ is an elementary dynamical piece. Moreover, $A$ is hyperbolic when $\zeta^\epsilon_p\in N_1$, or parabolic when $\zeta^\epsilon_p\in N_3$, with $\{\alpha(A),\omega(A)\}=\{\emptyset,\zeta^\epsilon_p\}$. 
	\item Suppose that $A$ is of type (2.d); that is, $A=\Sat_{F}(B)\cap U$, where $B$ is an open non-trivial interval in $\partial T_q$ with extremities $c_1,c_2\in\Upsilon_q\cup\Sing(\partial T_q)$, and $q\in N_3$. Moreover, since $B$ do not contain any singular point of the boundary $\partial T_q$, it is an analytic smooth curve. Assume that $q\in N_3^s$, for instance. Being $B$ contained in the frontier of $W^{tr}(p)=W^s_q$, the three-dimensional stable manifold at $q$, the positive leaf $\ell^+_x$ for any point $x\in\partial T_q$ (of type e-i relatively to $F$) is contained in $U$, except for the point $x$, and satisfies $\omega(\ell^+_x)=q$. Thus $A$ is saturated in $U$, an analytic submanifold,  and a parabolic elementary dynamical piece, with $\{\alpha(A),\omega(A)\}=\{\emptyset,q\}$. On the other hand, it is clear that $cl_U(A)\setminus A=\ell^U_{c_1}\cup\ell^U_{c_2}\cup\{q\}$, where $\ell^U_{c_i}:=\ell_{c_i}\cap U$. Taking into account that $c_1,c_2\in\Upsilon_q\cup\Sing(\partial T_q)$, we get that $\ell^U_{c_1},\ell^U_{c_2}$ are $U$-leaves and, in fact, 1-dimensional strata of types (1.c), (1.d) or (1.f). Thus, $cl_U(A)\setminus A$ is a union of strata of smaller dimension.   
	\item Finally, suppose that $\dim(A)=3$. By its very definition in (d3), $A$ is open (hence an analytic submanifold) and saturated in $U$, since $\mathcal{T}^2$ is saturated in $U$. The frontier $Fr_U(A)=cl_U(A)\setminus A$ is also saturated and contained in $\mathcal{T}^2$, thus formed of a union of strata of smaller dimension. Let us prove that $A$ is an elementary dynamical piece. Take one point $x_0\in A$. By Proposition~\ref{pro:omega-limits}, the $\omega$-limit set $\omega_U(\ell^+_{x_0})$ is either empty or a given singular point $q_0$. It suffices to prove that the set
	 $$B:=\{x\in A\,:\,\omega(\ell_x^+)=\omega(\ell^+_{x_0})\}$$
	  is equal to $A$ (we conclude analogously for $\alpha$-limit sets by taking $-\LL$ instead of $\LL$). Let us show that $B$ is open and closed in $A$. Assume first that $\omega(\ell^+_{x_0})=\emptyset$. Take $x\in B$. Then also $\omega_U(\ell^+_x)=\emptyset$. Since $x$ is an interior point of $U$, the first point, say $y$, in the positive leaf $\ell^+_{x}$ that belongs to $Fr(F)$ is of type i-t or i-e, relatively to $F$. The first case does not occur since those points belong to the family of intervals $\{L^i_p\}$ in Definition~\ref{def:fitting}, which would imply that a piece of the negative leaf through $y$ is included in $\mathcal{T}^2$, which is impossible. We get that $y$ belongs to the interior of some of the discs $T_p$, for some $p\in S_{tr}\cup N$, which permits to show that the set $B$ is open in $A$. On the other hand, suppose that $z\in A$ is the limit point of some sequence $\{x_n\}$ of points in $B$. Up to taking a subsequence, and by the above reasoning, each positive leaf $\ell^+_{x_n}$ cuts $Fr(F)$ for a first time at some point $y_n\in T_p$ for some $p\in S_{tr}\cup N$, independent of $n$. Being $T_p$ a compact set, we get that $\ell^+_z$ also cuts $T_p$, so that $\omega_U(\ell^+_z)=\emptyset$, and thus $z\in B$. This shows that $B$ is also closed in $A$. Finally, assume that $\omega(\ell^+_{x_0})=q_0$. Necessarily, we must have that $q_0\in N_3$ (and in fact in $N^s_3$): otherwise, the stable manifold at $q$, which contains the leaf $\ell^+_{x_0}$ locally at $q$, would be of dimension at most $2$ and contained in $\mathcal{T}^2$, which gives a contradiction. Thus $W^{tr}(q_0)=W^s_{q_0}$ is open and we obtain: if $x\in B$, then the positive leaf $\ell^+_x$ is either contained in the interior of $W^{tr}(q_0)$ (an open set) or cuts its boundary $\partial W^{tr}(q_0)$. We deduce similarly as above that $B$ is open and closed in $A$. 
\end{itemize}

\strut

Summarizing and coming back to the notation $U=U_n=int(F_n)$, where $F_n$ runs in the neighborhood basis of fitting domains, we have obtained a sequence $\{\SC_n\}_{n}$ of analytic $\LL$-stratifications of the neighborhood basis $\{U_n\}_n$.
From the description of the nature and closure of the different strata in the items above, statements (i) and (ii) in Theorem~\ref{th:afterblowup-intro} are proved straightforwardly.  On the other hand, taking into account that the number and properties of fixed and mobile separating surfaces (cf. Proposition~\ref{pro:fixed-separating} and Proposition~\ref{pro:mobile-separating}) do not depend on the fitting domain, we deduce that the number, nature and accumulation to the divisor of the strata listed in steps (d0)-(d3) do not depend on $n$, with the possible exception of those in (1.f) or (2.d). Since strata from (1.f) or (2.d) are parabolic, item (iii) is proved.

\vspace{.1cm}

This completes the proof of Theorem~\ref{th:afterblowup-intro}. \hfill{$\square$}

\strut

Now, it is clear how to prove the claims stated in Remark~\ref{rk:SC-topologic}. Notably, consider the stratification $\SC'_n$ of $U_n$, constructed as the stratification $\SC_n$, but with the following two modifications: 
\begin{itemize}
\item We get rid of those strata listed in item (1.f).
\item Consider a new family (2.d)' instead of the one in  (2.d) as follows: the intervals $B$ to be considered to produce the strata $\Sat_F(B)\cap U$ in the family (2.d)' are the connected components of $\partial T_q\setminus\Upsilon_q$, instead of those of $\partial T_q\setminus(\Upsilon_q\cup\Sing(\partial T_q))$. Take into account that $U$, the transversal disc $T_q$ at $q\in N_3$, and the set $\Upsilon_q$ depend on $n$, even if we have dropped the index $n$ to simplify the notation. However, the cardinality of $\Upsilon_q$ does not depend on the fitting domain, by property (v) in Definition~\ref{def:fitting}.
\end{itemize}
Then $\SC'_n=\{{A'}^n_j\}_{j\in J'}$ is an $\LL$-stratification of $U_n$ of class $\CC^0$ for which the set of indices $J'$ may be chosen independent of $n$ and satisfying the properties stated in Remark~\ref{rk:SC-topologic}: for each $j\in J'$, the character of ${A'}^n_j$, and the sets $\alpha({A'}^n_j),\omega({A'}^n_j),\overline{{A'}^n_j}\cap D$ do not depend on $n$. Besides, each strata of $\SC_n'$ is locally semi-analytic (i.e., semi-analytic at each of its points). Those strata given by (2.d)' (parabolic and two-dimensional) are the only ones which are possibly not smooth analytic submanifolds in the respective open set $U_n$. In fact, the analytic stratification $\SC_n$  described in Theorem~\ref{th:afterblowup-intro} is the refinement of $\SC'_n$ after removing the singular locus from  the non-smooth strata, that is, the one-dimensional strata in step (1.f) (see again Figure~\ref{Fig:RefinementStrata}). 

\begin{remark}{\em
The difference between the topological stratification $\SC'_n$ and its analytic refinement $\SC_n$ lies on the set of singular points in $\partial T_q$, with $q\in _3$. In turn, this singular set  strongly depends on the fitting domain choice. However, if the neighborhood basis $\{F_n\}_n$ can be built so that the number of such singular points may be bounded uniformly with respect to $n$, we can provide a sequence of analytic $\LL$-stratifications $\{\SC_n\}_n$ satisfying the stronger properties of the sequence $\{\SC_n'\}_n$ (namely, that item (iii) of Theorem~\ref{th:afterblowup-intro} is true for any stratum, including the parabolic ones).

If we recall the construction of the fitting domains $F_n$ made in \cite{Alo-San1}, we can see  that such a property holds: the fitting domain is obtained after a finite number of refinements from an initial ``distinguished fattening''  and, in each refinement, one produces at most a given, but controlled, number of new singular points in the frontier. We guess that one can even construct a sequence of fitting domains where the corresponding transversal discs have smooth bondary, but we have not completely checked this.}
\end{remark}

\section{Miscellaneous examples}\label{sec:miscellaneous}

\subsection{Some toy examples}
\begin{figure}[h]
\begin{center}
	\includegraphics[scale=0.6]{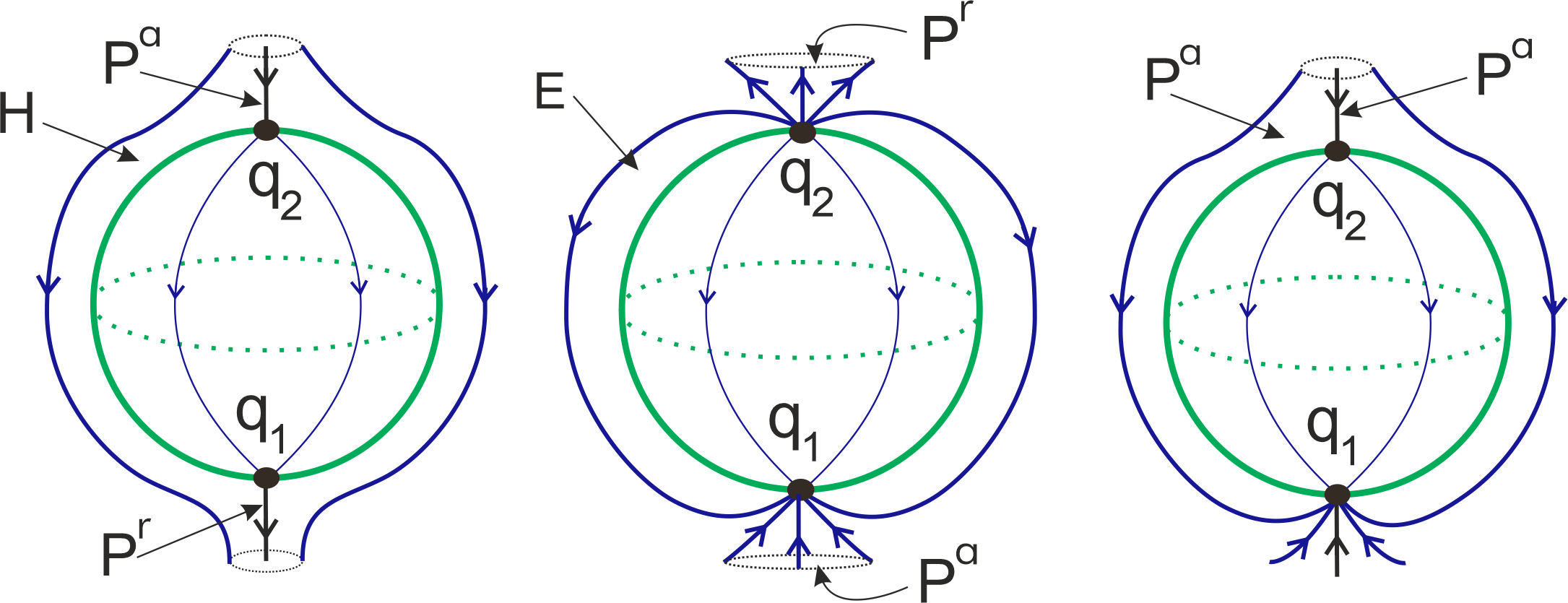}
	\caption{Possible strata configuration for two singular points on $D$. The notation $E,H,P^r,P^a$ means that the corresponding stratum is elliptic, hyperbolic, repelling parabolic and attracting parabolic, respectively.}
	\label{Fig:ToyExamples}
\end{center}
\end{figure}

The easiest situation in which we can provide the  stratification stated in Theorem~\ref{th:afterblowup-intro} corresponds to the case where the graph $\Omega$ of $\MM=(M,D,\LL)$ is just reduced to two isolated vertices $q_1,q_2$. According to \cite[Lemme 4]{Alo-San1}, except for a sign change, $q_1$ must be a repelling $D$-node and $q_2$ is an attracting one. Since there are no transversal saddle points, there are no separant surfaces and the stratification is easy to describe. In restriction to the divisor $D$ it is always given by the strata $\{q_1,q_2,D\setminus\{q_1,q_2\}\}$, with the last one being elliptic.  However, outside $D$ the stratification depends on  the nature of the points $q_1,q_2$, that is, if they are three-dimensional saddles or not. In Figure \ref{Fig:ToyExamples} we depict the three situations: 

(a) $q_1,q_2\in N_1$. 

(b) $q_1,q_2\in N_3$. 

(c) $q_1\in N_3$, $q_2\in N_1$  (analogously $q_1\in N_1$, $q_2\in N_3$).

At this point, we can ask if a foliated manifold $\MM=(M,D,\LL)$ in the above situations could be realized from germs of analytic vector fields $\xi$ at $0\in\R^3$ after performing the polar blowing-up
$$
\pi_0:\R_{\ge 0}\times\SSS^2\to\R^3,\;\pi(r,X)=rX.
$$
That is, we would have $M=\R_{\ge2}\times\SSS^2$, $D=\{0\}\times\SSS^2$ and $\LL$ generated by the strict transform $r^{-(\nu-1)}\pi^*\xi$, where $\nu$ is the multiplicity of $\xi$ at $0$.

Examples for situations (a) and (b) are easy to obtain. For instance, we can take the following vector fields, respectively:
$$
\begin{array}{l}
\xi_a:=-(x^2+y^2+z^2)\frac{\partial}{\partial z}
\\
\xi_b:=2xz\frac{\partial}{\partial x}
+
2yz\frac{\partial}{\partial y}
+
\left(z^2-(x^2+y^2)\right)\frac{\partial}{\partial z}
\end{array}
$$
In both cases, $q_2$ (the repelling $D$-node) is the north pole of the sphere, whereas $q_1$ (the attracting $D$-node) is the south pole. 

We have not found an example of a vector field at $0\in\R^3$ giving rise to situation (c) after blowing-up. To be honest, we are not sure whether it exists or not.

\subsection{The graph $\Omega$ and the  transversal eigenvalues signs do not determine the stratification}

Now we propose an example that permits to show that, contrary to what happens in dimension two (see Section \ref{sec:intro}), the nature of the strata as well as  their adjacency relations are not univocally determined by the graph $\Omega$ and the character stable/unstable of the transversal manifolds $W^{tr}(p)$ at the singular points $p\in S'=N\cup S_{tr}$.

\begin{figure}[h]
	\begin{center}
		\includegraphics[scale=0.6]{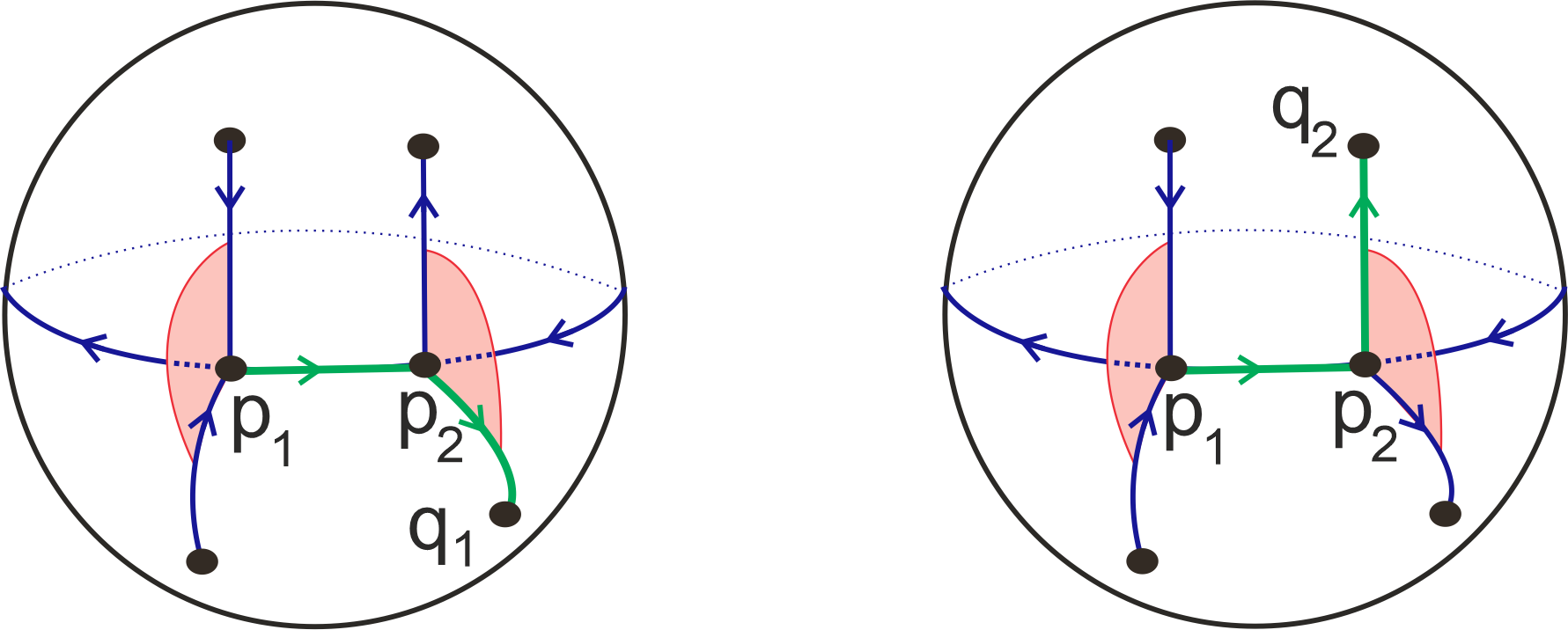}
		\caption{Different strata configuration with identical graph for $p_1\in S_{tr}^s$, $p_2\in S_{tr}^u$, $q_1\in N_3^s$, $q_2\in N_1^u$.}
		\label{Fig:ToyExamples2}
	\end{center}
\end{figure}

Consider, for instance, the two situations depicted in Figure~\ref{Fig:ToyExamples2}. In both cases, the graph $\Omega$ and the partition $V(\Omega)=N_1\cup N_3\cup S_{tr}^s\cup S_{tr}^u$ of the set of vertices (singular locus) are identical. However, the traces along the divisor of mobile separating surfaces generated at one of the transversal saddle points, say $p_1$, are different in each of the situations.  Such trace is the path of edges $\Theta(\nu_{p_1})$ described in Proposition~\ref{pro:saturation-marks} for one of the local sides $\nu_{p_1}$ of $W^{tr}(p_1)$. Notably, in the first case, $\Theta(\nu_{p_1})$ ends at the $D$-node point $q_1\in N_3$ while, in the second case, it ends at $q_2\in N_1$. Taking into account the strata description in Section~\ref{sec:proof}, we deduce that the open stratum lying ``between'' the two corresponding mobile separating surfaces is parabolic attracting in the first situation but hyperbolic in the second one.

The conclusion is that the knowledge of the traces of mobile and separating surfaces is crucial to describe the strata configuration. These traces depend, not only on the signs of transversal eigenvalues, but also in the verification that the structural weight assignments at singular points  satisfy, or not, a certain system of polynomial inequalities (cf. discussion of s-resonances in Section~\ref{sec:background}). 

It is worth to remark that, for these examples,  we have not tackled the question of whether they are realizable from a local analytic vector field or even as a foliated variety $(M,D,\LL)$ with $M=\R_{\ge0}\times\SSS^2$. Strictly speaking, our example is thus hypothetical, but representative of a phenomena that, certainly, seems to occur in true examples. 

\subsection{The germs of hyperbolic strata are not well defined}
As announced in Section \ref{sec:intro}, we propose here an example of a vector field at $0\in\R^3$ which is stratifiable on each element of a basis of neighborhoods $\{U_n\}$ such that the germs of the hyperbolic strata of the different elements $U_n$ do not coincide. More properly, if $\HH_n\subset U_n$ denotes the union of the hyperbolic strata in $U_n$ for each $n$, then the equation $\HH_n\cap U_{n'}=\HH_{n'}\cap U_n$ do not hold necessarily for $n\ne n'$.

\begin{figure}[h]
\begin{center}
	\includegraphics[scale=0.7]{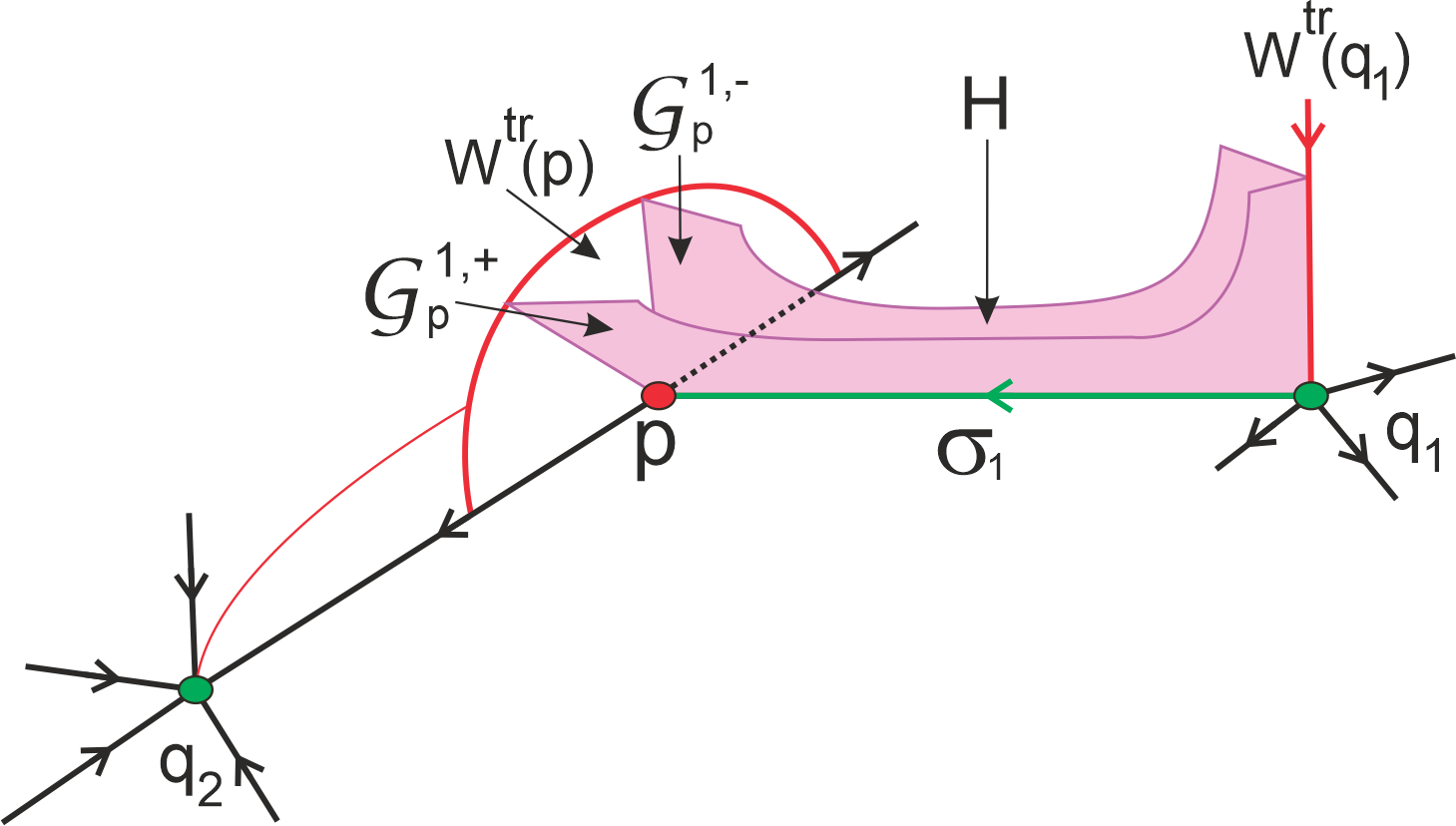}
	\caption{Configuration $C$.}
	\label{Fig:ConfigurationC}
\end{center}
\end{figure}

Setting the same problem in the general context of foliated manifolds $\MM=(M,D,\LL)$, we obtain such an example (when the hypothesis in Theorem~\ref{th:afterblowup-intro} hold) if the associated graph $\Omega_\MM$ contains the following configuration (see Figure~\ref{Fig:ConfigurationC}):
\begin{quote}
\!\!\! (C)\;\; In the interior of one component of $D$ there are singular points $q_1\in N_1$, $p\in S_{tr}$ and $q_2\in N_3$ and edges $\sigma_1$, $\sigma_2$ with extremities $q_1,p$ and $q_2,p$, respectively, such that $\sigma_1$ contains the $1$-dimensional invariant manifold $W^1_p$ at $p$, whereas $\sigma_2$ contains $W^{tr}(p)\cap D$. 
\end{quote}

\noindent Notice that if we have the situation (C) and $F$ is a fitting domain of $\MM$, the two mobile separating surfaces $\GG_p^{1,-},\GG_p^{1,+}$ generated at $p$ and associated to the local side of $W^{tr}(p)$ that contains $\sigma_1$, enclose an open hyperbolic stratum $H$ (composed of leaves that escape from $F$ along the transversal discs $T_{q_1}$ and $T_p$ of the fitting domain). Also, one of those separating surfaces, say $\GG_p^{1,+}$, is contained in the boundary of a   
parabolic stratum $P$ (composed of leaves escaping $F$ in one sense through $T_{q_1}$ and accumulating to $q_2$ in the other one). But the mobile separating surfaces $\GG_p^{1,-},\GG_p^{1,+}$ (hence $H$ and $P$) depend strongly on the fitting domain $F$. In particular, those surfaces for  the different $F$ inside a neighborhood basis may not ``germify'' along the divisor $D$.

\vspace{0.2cm}
Consider the following polynomial vector field in $\R^3$, depending on a parameter $\alpha\in\R$:
\begin{equation}\label{eq:example-hyperbolic}
\xi_\alpha:=x^2(z-\alpha y+\frac{1}{2}x)
\frac{\partial}{\partial x}+
y(xz-\frac{1}{2}x^2+y^3)
\frac{\partial}{\partial y}+
z(-z^2+zx-yx+\frac{3}{2}x^2)
\frac{\partial}{\partial z}.
\end{equation}

{\bf Claim.-} If $\alpha<1$ sufficiently closed to $1$,  there exists a reduction of singularities $\pi:M\to\R^3$ of $\xi_\alpha$, consisting of finitely many punctual blowing-ups such that the 
foliated manifold $\MM=(M,\pi^{-1}(0),\pi^*\LL_{\xi_\alpha})$ is hyperbolic, acyclic, non-dicritical, of Morse-Smale type and with no s-resonances and such that the correspnding graph $\Omega_{\MM}$ contains a configuration of the form (C).

\strut

\begin{figure}[h]
	\begin{center}
		\includegraphics[scale=0.6]{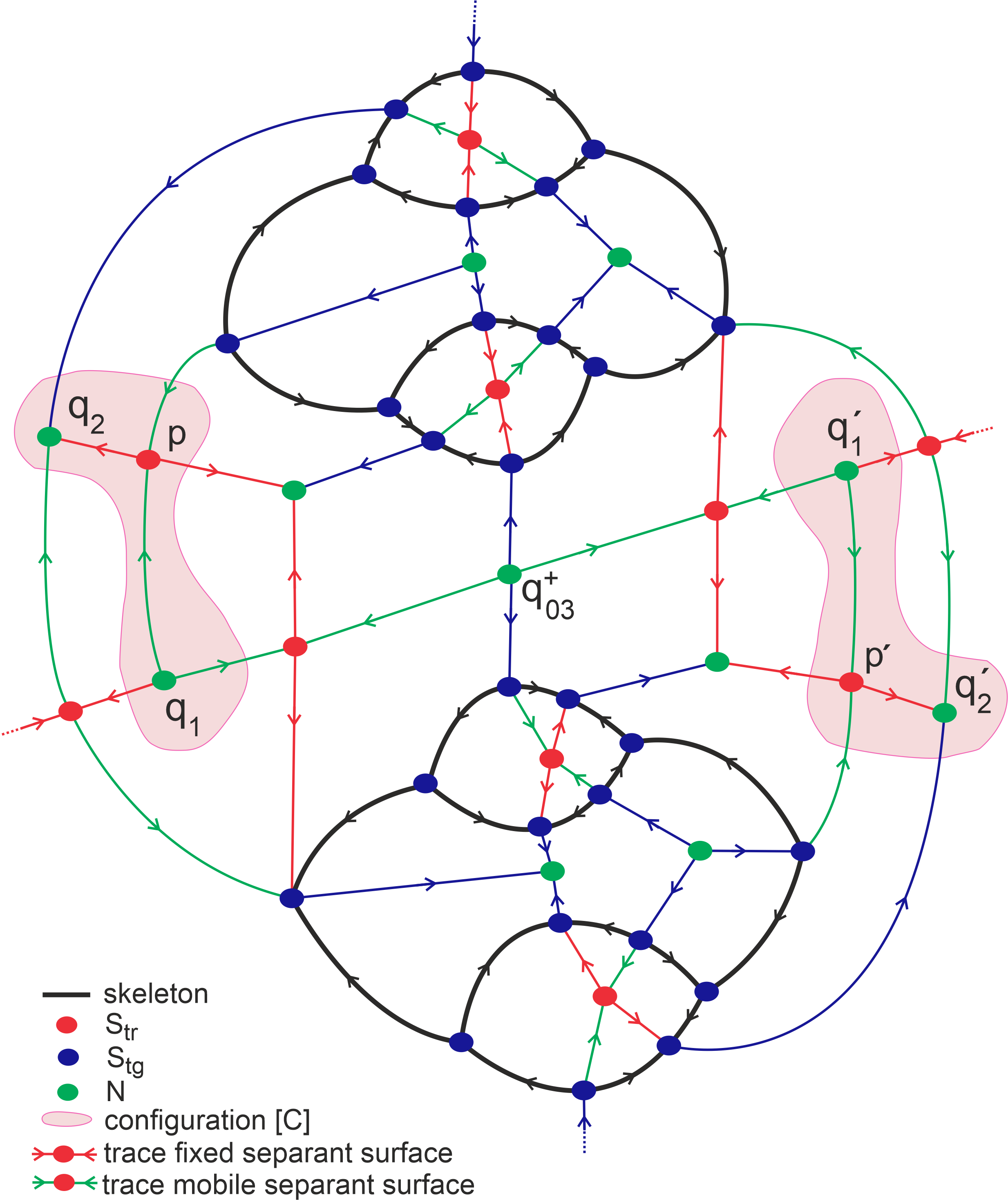}
		\caption{The graph $\Omega$ corresponding to the last example.}
		\label{Fig:LastExample}
	\end{center}
\end{figure}

We do not reproduce here all the calculations needed to show the claim. We just indicate to the reader the sequence of blowing-ups giving rise to the morphism $\pi$.
\begin{enumerate}
	\item  We consider the polar blowing-up of the origin $\pi_0:M_0=\SSS^2\times\R_{\ge0}\to\R^3$. Instead of using spherical coordinates, we consider $M_0$ covered by the six charts representing the six half-planes associated to the coordinates planes $(x=0)$, $(y=0)$, $(z=0)$. Namely, we have charts 
	$$
	\left\{
	(U^\eps_x,(x_1^\eps,y_1^\eps,z_1^\eps)),
	(U^\eps_y,(x_2^\eps,y_2^\eps,z_2^\eps)),
	(U^\eps_z,(x_3^\eps,y_3^\eps,z_3^\eps))
	\right\}_{\eps\in\{+,-\}},
	$$
	where the range of the variables is: $\R_{\ge0}$ for $x_1^\eps,y_2^\eps,z^\eps_3$ and $\R$ for the rest. The respective expressions of the map $\pi_0$, together with the equation of the divisor $D_0=\pi^{-1}(0)$, are given by
	$$
	\begin{array}{c}
		\pi_0(x_1^\eps,y_1^\eps,z_1^\eps)=(\eps x_1^\eps,x_1^\eps y_1^\eps,x^\eps_1z_1^\eps),\;\;D_0\cap U_x^\eps=(x_1^\eps=0)\\
		\pi_0(x_2^\eps,y_2^\eps,z_2^\eps)=(\eps y_2^\eps x_2^\eps,\eps y_2^\eps,y_2^\eps z_2^\eps),\;\;D_0\cap U_y^\eps=(y_2^\eps=0)\\
		\pi_0(x_3^\eps,y_3^\eps,z_3^\eps)=(z_3^\eps x_3^\eps,z_3^\eps y_3^\eps,\eps z_3^\eps),\;\;D_0\cap U_z^\eps=(z_3^\eps=0).
	\end{array}
	$$ 
	\item We compute the strict transforms of $\xi_\alpha$ in each chart (i.e.,  $\frac{1}{(x_1^\eps)^2}(\pi_0|_{U_1^{\eps}})^*\xi_\alpha$, etc.) and we obtain finitely many singular points of the transformed foliation $\pi_0^*\LL_{\xi_\alpha}$ on $D_0$. More precisely:
	\begin{itemize}
		\item In each one of the charts $U^+_x,U^-_x$, we get six singular points, all of them hyperbolic, being three of them  in the configuration (C). This can be checked more easily for $\alpha=1$ and conclude that it is also true for a perturbation $\alpha\approx 1$. 
		\item In $U^\eps_z\setminus(U^+_x\cup U^-_x)$, there is only one singular point, the origin $q_{03}^\eps$ of the chart $U^\eps_z$, which is hyperbolic (in fact, it belongs to $N_1$).
		\item The origins $q_{02}^\eps$ of the charts $U^\eps_y$ (the only points not covered by the other charts) are also singular, but not hyperbolic.
	\end{itemize}
\item  We perform new polar blowing-ups at the points $q^+_{02}$ and $q^-_{02}$. For instance, for the blowing-up $\pi^+_1$ at $q^+_{02}$, we consider new charts $U_{x_2}^\eps$, $U^+_{y_2}$, $U^{\eps}_{z_2}$ with respective coordinates $(x^\eps_{11},y^{\eps}_{11},z^{\eps}_{11})$, etc. (the chart $U^-_{y_2}$ is not longer valid since the coordinate $y^+_2$ at $q^{+}_{02}$ takes only values in $\R_{\ge0}$). After computing the corresponding strict transforms in the different charts, we get that the transformed foliation $(\pi_0\circ\pi^+_1)^*\LL_{\xi_\alpha}$ has only finitely many singular points along the new component of the divisor $D_1^+:=(\pi^+_1)^{-1}(q^+_{02})$ (here, the condition that $\alpha\ne 1$ is necessary). Moreover, they are hyperbolic except for the origins $q^\eps_{13}$ of the charts $U^{\eps}_{z_2}$. Similar conclusions are obtained when we perform the blowing-up $\pi_1^-$ at the point $q^-_{02}$.

\item  We have to continue making new blowing-ups $\pi^\eps_2$ at the points $q^\eps_{13}$ (and the corresponding ones coming out from $\pi^-_1$). After such blowing-ups, we finally get a reduction of singularities $\pi:M\to\R^3$ of $\xi_\alpha$, since there are only finitely many singular points along the new components $D^\eps_{2}:=(\pi^\eps_2)^{-1}(q^\eps_{13}))$. All of them are hyperbolic.

\item Computing the signs of all hyperbolic singular points appearing after the above sequence of blowing-ups $\pi$, we obtain the graph $\Omega=\Omega_\MM$ associated to the corresponding foliated manifold $\MM=(M,D=\pi^{-1}(0),\pi^*\LL_{\xi_\alpha})$. It is despicted in Figure~\ref{Fig:LastExample}.  In that picture, we have distinguished the singular points according to whether they are transveral saddles, tangential saddles or $D$-nodes (we have not specified if the latter ones belong either to $N_1$ or to $N_3$). Also, we have indicated the skeleton of the manifold $M$, as well as the traces of the different fixed and mobile separating surfaces. To make it easier to read, we have simplified the picture sending to infinity the origin $q^-_{03}$ of the chart $U^-_z$ of the first blowing-up $\pi_0$. As mentioned before, this point is a $D$-node point (in fact, it belongs to $N_1^{s}$). Thus, the boundedless edges of the picture represent edges starting at this point. 

One can check that the foliated manifold $\MM$ is acyclic provided that the graph $\Omega$ has not polycycles and, since any cycle must contain a singularity inside, it would cut  $|\Omega|$ in case it exists.
	
It remains to show that  $\MM$ is also of Morse-Smale type and has no s-resonances. The first property, being of Morse-Smale type, is direct from Figure~\ref{Fig:LastExample} (we only need to know the graph and what edges are included in the skeleton). Concerning the second property, no s-resonances, despite in general it does not follow from the graph alone $\Omega$ (we need to check certain algebraic inequalities for the eigenvalues at singular points), in this particular example it does. We can check from Figure~\ref{Fig:LastExample} that we do not  have multiple saddle connections $(\sigma_0,\dots,\sigma_n)$ along the skeleton such that the first and last edges $\sigma_0,\sigma_n$ support the one-dimensional invariant manifolds of the corresponding singular points (cf. definition of s-resonance in Subsection ~\ref{sec:additional}).
	
\end{enumerate}

\end{document}